\documentclass[12pt, twoside, leqno]{article}

\usepackage{amsmath,amsthm}
\usepackage{amssymb}

\usepackage{enumitem}

\usepackage{graphicx}

\usepackage[T1]{fontenc}

\newtheorem{theorem}{Theorem}[section]

\newtheorem{lemma}[theorem]{Lemma}

\theoremstyle{definition}

\numberwithin{equation}{section}

\usepackage{amsmath,amssymb,amsthm}
\usepackage{mathrsfs}    
\usepackage{mathtools}   
\usepackage{yfonts}      
\usepackage{dirtytalk}   
\usepackage{tikz-cd}     
\usepackage{hyperref}    

\theoremstyle{plain}
\newtheorem*{theoremA*}{Theorem A}
\newtheorem*{theoremB*}{Theorem B}
\newtheorem*{theoremC*}{Theorem C}
\newtheorem*{theoremD*}{Theorem D}
\newtheorem*{theorem*}{Theorem}

\begin{document}

\title{Eigenspaces of Coleman's Trace Operator}


\author{Joseph DiCapua\\
Department of Mathematics, Louisiana State University\\ 
175 Field House Drive, Baton Rouge 70802, LA, United States\\
E-mail: jdicap1@lsu.edu}

\date{}

\maketitle


\renewcommand{\thefootnote}{}

\footnote{2020 \emph{Mathematics Subject Classification}: Primary 11S99.}

\footnote{\emph{Key words and phrases}: Lubin--Tate formal group laws, Iwasawa theory, Coleman power series, dihedral $(\varphi,\Gamma)$-modules, Iwasawa cohomology.}

\renewcommand{\thefootnote}{\arabic{footnote}}
\setcounter{footnote}{0}


\begin{abstract}
The Coleman power series defined on a Lubin--Tate tower of extensions over $K$ are compatible with respect to two formal group laws: the multiplicative formal group law and some Lubin--Tate formal group law defined over $\mathcal{O}_K$. We ask if it is possible to generalize Coleman series in order to find series which are compatible with respect to any two Lubin--Tate formal group laws in the same way. We provide a precise formulation of this question and a partial answer towards the classification of all such power series that involves the eigenspaces of Coleman's trace operator. Some additional eigenspaces of Coleman's trace operator are also introduced. These new eigenspaces are used to construct certain fixed points of the $\psi$-operator on the family of trianguline dihedral $(\varphi, \Gamma)$-modules.
\end{abstract}

\section{Introduction}

\noindent
Understanding ideal class groups of number fields is one of the most important problems in number theory. Iwasawa theory originated when Kenkichi Iwasawa investigated ideal class groups in cyclotomic towers and their relationships to $p$-adic $L$-functions starting in the late 1950s. These relationships are known as \say{Iwasawa main conjectures.} Since the original \say{Iwasawa main conjecture} was proved in the 1980s by Mazur--Wiles, there has been significant interest in generalizing these theorems to different settings (elliptic curves, modular forms, etc.) \cite{Ochiai}. Many of these generalizations use Euler systems. Such constructions have given deep results towards understanding Tate--Shafarevich groups and the Birch and Swinnerton-Dyer conjecture.
\\
\\
Iwasawa theory also has important connections to the local Langlands conjectures. One way this occurs is through the study of $(\varphi,\Gamma)$-modules, and the $\psi$-operator on these modules. The fixed points of this operator are closely related to Iwasawa cohomology and $p$-adic $L$-functions. For the $p$-adic local Langlands correspondence constructed by Pierre Colmez in \cite{ColmezRepDeGL2}, \say{The classical local Langlands correspondence is normalized via the local factors of $L$-functions, whereas the $p$-adic correspondence, described above, relies on a direct construction. However, the local factors of $p$-adic $L$-functions play a very important hidden role.} He states later in the same paper: \say{One can therefore hope that the [fixed points of the $\psi$-operator] will play an important role in the $p$-adic local Langlands correspondence. This is indeed the case, and ideas from Iwasawa theory intervene at several places in the determination of locally algebraic vectors.}
\\
\\
For any finite Galois extension of $p$-adic number fields $M/K$, one has the classical norm map $\text{N}_{M/K}(x)$ defined by multiplying the conjugates of $x$. For any Lubin--Tate formal group law $F$, one can consider an analog of this map where we replace multiplication with the operation coming from the formal group law $F$. When $F$ is a Lubin--Tate formal group law associated to a proper subfield of $K$, we call such analogs $f_L$-norms. An \hyperref[alphaf_0normcompatibledef]{$\alpha\text{-}f_L\text{-norm}$ compatible sequence} is a sequence such that consecutive entries are compatible with respect to some $f_L$-norm with coefficient $\alpha$. The first goal of this paper is to classify certain power series which interpolate $\alpha\text{-}f_L\text{-norm}$ compatible sequences. We are able to give an answer towards classifying all such power series in \hyperref[sec:section2]{section 2} using certain eigenspaces of Coleman's trace operator. The eigenvalues for these eigenspaces are all constants. 
\\
\\
The second goal of this paper is to introduce a new class of eigenspaces of Coleman's trace operator, or the $\psi$-operator from $p$-adic Hodge theory, for which the eigenvalues are series, rather than constants, and apply these eigenspaces to the study of Iwasawa cohomology of trianguline dihedral $p$-adic Galois representations. We introduce new eigenspaces in \hyperref[sec:section3]{section 3} and give an application to trianguline dihedral $(\varphi,\Gamma)$-modules in \hyperref[sec:section4]{section 4}. In particular, given a module in a certain family of trianguline dihedral $(\varphi,\Gamma)$-modules, we are able to construct all fixed points of the $\psi$-operator satisfying particular conditions.
\\
\\
Robert Coleman's work \say{Division Values in Local Fields} \cite{Coleman1} introduced Coleman power series which are ubiquitous in Iwasawa theory. Some applications of Coleman power series include the Coates-Wiles homomorphisms and more general Coleman maps, which are used to construct $p$-adic $L$-functions. They are also used in one proof of the classical main conjecture of Iwasawa theory. Since Coleman power series are of great importance in Iwasawa theory, it is a natural question to ask for generalizations of Coleman's construction. This paper provides a family of modules of power series which generalize Coleman power series in that they are compatible with respect to two Lubin--Tate formal group laws. The aforementioned modules are also closely related to eigenspaces of Coleman's trace operator, and we show eigenspaces of this operator can be used to construct fixed points of the $\psi$-operator for certain $(\varphi,\Gamma)$-modules.
\\
\\
It should be noted that the questions considered in this paper are also motivated by the question asked by Laurent Berger and Lionel Fourquaux in \cite[Remark 3.4.7]{BergerFourquaux} and its solution in \cite{DiCapuaKolyvagin}. The series the author and Victor Kolyvagin classified in \cite{DiCapuaKolyvagin} generalize Coleman series in the cyclotomic case to the Lubin--Tate case. However, Coleman power series are more general, as they allow one to parametrize all norm compatible sequences over any tower of Lubin--Tate extensions. In this paper we develop a further generalization of Coleman series in which for any pair of Lubin--Tate formal group laws (one defined over a base field $L$ and another defined over $K$, an extension of $L$), we describe a family of modules of power series compatible with respect to both of these formal group laws.
\\
\\
We introduce a \hyperref[eq:alphamodule]{family of $\mathcal{O}_L$-modules $\mathscr{A}^\alpha_{L,K}$}, which consist of power series interpolating $\alpha\text{-}f_L\text{-}$norm compatible sequences. Here $K$ is a finite extension of $L$, and $\alpha$ is a parameter defining the recursive relation satisfied by the sequences. When $L = \mathbb{Q}_p$ and $\alpha = 1$, the module $\mathscr{A}^\alpha_{L,K}$ is isomorphic to the $\mathbb{Z}_p$-module of Coleman power series interpolating norm compatible sequences of principal units defined over $K$. When $L = K$ and $\alpha = q_L/\pi_L$ the module $\mathscr{A}^\alpha_{L,K}$ is isomorphic to the module $\mathscr{A}$ of all power series $f(x) \in \mathcal{O}_K[[x]]$ interpolating Lubin--Tate trace compatible sequences from \cite{DiCapuaKolyvagin,BergerFourquaux}. For different values of $\alpha$ and different choices of fields $L$ and $K$ the modules $\mathscr{A}^\alpha_{L,K}$ are analogs of Coleman power series. A central result of this paper is to give an answer towards classifying the power series in $\mathscr{A}^\alpha_{L,K}$.
\\
\\
We are able to prove the following result which is analogous to the main result from \cite{DiCapuaKolyvagin}:

\begin{theoremA*}
    (See \hyperref[thm:2.0.2]{Theorem 2.2}) For values $\alpha$ dividing $q_L/\pi_L$ there is a map from $\mathscr{A}^\alpha_{L,K}$ to the kernel of Coleman's trace operator. The kernel of this map is either trivial or generated by a single series as an $\mathcal{O}_L$-module. For all but one such value of $\alpha$ the above map is injective, so we conclude that $\mathscr{A}^\alpha_{L,K}$ is isomorphic to some submodule of the kernel of Coleman's trace operator in these cases.
\end{theoremA*}

\noindent
We also make use of eigenspaces of Coleman's trace operator to give a partial answer towards classifying all series in the modules $\mathscr{A}^\alpha_{L,K}$. Let $\pi_L$ be a fixed uniformizer of $\mathcal{O}_L$. Let $\mathscr{A}^\alpha_{L,K,\pi_L}$ denote the submodule of $\mathscr{A}^\alpha_{L,K}$ consisting of all $r(x) \in \mathscr{A}^\alpha_{L,K}$ such that $r(x) \in \pi_L\mathcal{O}_K[[x]]$. Denote the $\alpha$-eigenspace of Coleman's trace operator in $\pi_L\mathcal{O}_K[[x]]$ by $\mathscr{E}^\alpha_{L,K}$. We have the following:

\begin{theoremB*}
    (See \hyperref[thm:2.0.3]{Theorem 2.3}) Let $\log_{F_L}$ be the logarithm of the formal group law $F_L$. Then the map $\log_{F_L} : \mathscr{A}^\alpha_{L,K,\pi_L} \rightarrow \mathscr{E}^\alpha_{L,K}$ is an isomorphism of $\mathcal{O}_L$-modules.
\end{theoremB*}

\noindent
The $\alpha$-eigenspace of Coleman's trace operator was solved for in \cite{DiCapuaKolyvagin}, so the above result gives a large collection of power series in $\mathscr{A}^\alpha_{L,K}$ when $\pi_L \mid \alpha$. We use $\mathscr{L}_{F_K}$ to denote Coleman's trace operator defined with respect to the Lubin--Tate formal group law $F_K$. We provide additional proofs that the $\alpha$-eigenspace of $\mathscr{L}_{F_K}$ is isomorphic to the kernel of $\mathscr{L}_{F_K}$.
\\
\\
In \hyperref[sec:section3]{section 3} of this paper we present some additional eigenspaces of Coleman's trace operator. If $\alpha \in \mathcal{O}_K$ then the $\alpha$-eigenspace of $\mathscr{L}_{F_K}$ has already been found as mentioned above. In \hyperref[sec:section3]{section 3} we show one can replace the constant $\alpha$ with a power series $\alpha(x) \in \pi_K^2\mathcal{O}_K[[x]]$. $\mathscr{E}^\alpha_K$ is defined to be the $\mathcal{O}_K$-module of all series $f(x) \in \mathcal{O}_K[[x]]$ such that
$$\mathscr{L}_{F_K}(f) = \alpha(x)f(x)$$

\noindent
We get the following:

\begin{theoremC*}
    (See \hyperref[thm:3.0.1]{Theorem 3.1}) There is an explicit map $\mathscr{E}^\alpha_K \rightarrow \mathscr{C}_K$, where $\mathscr{C}_K$ denotes the kernel of $\mathscr{L}_{F_K}$ in $\mathcal{O}_K[[x]]$. This map is an isomorphism of $\mathcal{O}_K$-modules.
\end{theoremC*}

\noindent
In \hyperref[sec:section4]{section 4} of this paper we give one application of the eigenspaces from \hyperref[sec:section3]{section 3} to the family of trianguline dihedral $(\varphi, \Gamma)$-modules. We study the description of the $(\varphi, \Gamma)$-modules associated to different dihedral representations from \cite{Berger1}. For the particular basis used by Berger in the same paper and under certain conditions, we are able to give a description of all fixed points of the $\psi$-operator which have coefficients in the integral part of the Robba ring. 
\\
\\
Let $\mathcal{R}$ denote the Robba ring over $\mathbb{Z}_p$, so $\mathcal{R}$ consists of series of the form $f(x) = \sum_{n \in \mathbb{Z}} a_n x^n$ such that $a_n \in \mathbb{Q}_p$ and $f$ converges on some annulus $r_f \le |x| <1$ depending on $f$. We define the Frobenius map $\varphi$ on $\mathcal{R}$ so that $\varphi(f)(x) = f((1+x)^p-1)$. Let $Q = \varphi(x)/x$, and let $Q_2 = \varphi(Q)$. We prove the following:

\begin{theoremD*}
    (See \hyperref[thm:4.0.4]{Theorem 4.4}) Let $D$ be a $(\varphi,\Gamma)$-module over $\mathcal{R}$ having a basis $\{e_1,e_2\}$ such that the matrix for $\varphi$ is given by

    $$\text{\emph{Mat}}(\varphi) = \begin{pmatrix}
        0 & 1\\
        Q_2^r(Q_2/Q^p)^u & 0
    \end{pmatrix}
    $$

    \noindent
    Then each fixed point of $\psi$ on $D$ with integral coefficients with respect to the basis $\{e_1,e_2\}$ corresponds to a series in the $Q_2^r(Q_2/Q^p)^u$-eigenspace of the $\psi^2$-operator on the integral part of $\mathcal{R}$. 
\end{theoremD*}

\noindent
We are then able to completely describe the $Q_2^r(Q_2/Q^p)^u$-eigenspace of the $\psi^2$-operator on $\mathcal{R}^\text{int}$, the integral part of $\mathcal{R}$, whenever $u \in \mathbb{Z}$ and $s = r +u(p+1) \le 0$. For the cases when $u \in \mathbb{Z}$ and $s > 0$ we are able to determine a submodule of the $Q_2^r(Q_2/Q^p)^u$-eigenspace of the $\psi^2$-operator on $\mathcal{R}^{\text{int}}$. Note that by \cite{Berger1}, the condition $r,u \in \mathbb{Z}$ is equivalent to the corresponding dihedral representation being trianguline, so these results apply to trianguline dihedral $p$-adic Galois representations.
\\
\\
The results from \hyperref[sec:section4]{section 4} raise several natural questions. These include whether the maps and results presented in that section can be used to parametrize classes in Iwasawa cohomology across the entire family of trianguline dihedral $p$-adic Galois representations in a useful way and whether the results from the same section generalize beyond the trianguline case.

\subsection{Notation}

Throughout this paper $L$ is a finite extension of $\mathbb{Q}_p$. $\mathcal{O}_L$ denotes the ring of integers in $L$. $\pi_L$ denotes a uniformizer for $\mathcal{O}_L$. $q_L = |\mathcal{O}_L/\pi_L\mathcal{O}_L|$ denotes the size of the residue field of $\mathcal{O}_L$. We fix a series $f_L \in \mathcal{O}_L[[x]]$ such that $f_L(x) \equiv x^{q_L} \mod \pi_L$ and such that $f_L(x) \equiv \pi_Lx \mod x^2$. We denote the Lubin--Tate formal group law associated to $f_L$ by $F_L$. If $a \in \mathcal{O}_L$, then the associated endomorphism of $F_L$ is denoted by $[a]_L(x)$. For $x$ and $y$ having positive valuation, or $x$ and $y$ being power series having constant term with positive valuation, we will write $x\oplus_L y = F_L(x,y)$ to denote the formal group law operation with respect to $F_L$.
\\
\\
Let $K$ be a finite extension of $L$. Fix a uniformizer $\pi_K$ of $\mathcal{O}_K$ and let $q_K$ denote the size of the residue field of $\mathcal{O}_K$. Fix a series $f_K \in \mathcal{O}_K[[x]]$ such that $f_K(x) \equiv x^{q_K} \mod \pi_K$ and such that $f_K(x) \equiv \pi_Kx \mod x^2$. We denote the Lubin--Tate formal group law associated to $f_K$ by $F_K$. Let $I_{\overline{K}}$ denote the maximal ideal of the ring of integers of the algebraic closure of $K$. Then $I_{\overline{K}}$ has the structure of an $\mathcal{O}_K$-module. We let $x\oplus_K y = F_K(x,y)$ whenever $x,y \in I_{\overline{K}}$. For any $a \in \mathcal{O}_K$ we denote the associated endomorphism of $F_K$ by $[a]_K(x) \in \mathcal{O}_K[[x]]$.
\\
\\
We fix a sequence $(u_n)_{n\ge 0}$ such that each $u_n \in I_{\overline{K}}$ is a root of $f_K^{(n+1)}(x)$. We also stipulate that for each $n \ge 0$ we have $[\pi_K]_K(u_{n+1}) = f_K(u_{n+1}) = u_n$ and $u_0 \neq 0$. These conditions guarantee that each $u_n$ generates the roots of $f_K^{(n+1)}(x)$ as an $\mathcal{O}_K$-module. We let $K_n = K(u_n)$. We let $z_1, z_2, \ldots, z_{q_K-1}$ be the nonzero roots of $f_K(x)$.
\\
\\
For any finite extension $M$ of $\mathbb{Q}_p$ we let $I_M$ denote the maximal ideal of $\mathcal{O}_M$, the ring of integers in $M$. For each $n \ge 0$ we define a map $f_L\text{-norm}_{K_{n+1}/K_n}: I_{K_{n+1}} \rightarrow I_{K_n}$ in the following way: for each $x \in I_{K_{n+1}}$ we let
$$f_L\text{-norm}_{K_{n+1}/K_n}(x) = \sigma_1(x) \oplus_L \sigma_2(x) \oplus_L \ldots \oplus_L \sigma_{q_K}(x)$$

\noindent
where $\{\sigma_i\} = \text{Gal}(K_{n+1}/K_n)$. We use the notation $f_L$-norm to emphasize that the Lubin--Tate formal group law addition being used is not attached in any way to the field extensions $K_n$. This is to distinguish our notation from the Lubin--Tate trace in \cite{BergerFourquaux}, which is primarily used in that paper for a Lubin--Tate formal group law attached to a uniformizer of $K$. 
\\
\\
Note that in the special case $K=L$ and $f_K = f_L$ we have $f_L\text{-norm}_{K_{n+1}/K_n}$ is the same map as the Lubin--Tate trace defined in \cite{BergerFourquaux}. We then fix some $\alpha \in \mathcal{O}_L$ and we define an $\alpha\text{-}f_L\text{-norm}$ compatible sequence to be any sequence $(x_n)$ \label{alphaf_0normcompatibledef} with each $x_n \in I_{K_n}$ such that 
$$f_L\text{-norm}_{K_{n+1}/K_n}(x_{n+1}) = [\alpha]_L(x_n)$$

\noindent
for each $n \ge 0$. We say that such a sequence is parametrizable or interpolated if there exists $r(x) \in \mathcal{O}_K[[x]]$ such that $r(u_n) = x_n$ for each $n \ge 0$.

\section{Interpolating $f_L$-norm compatible sequences}
\label{sec:section2}
\noindent
For different values of $\alpha \in \mathcal{O}_L$ one can ask whether it is possible to classify all interpolated $\alpha\text{-}f_L\text{-norm}$ compatible sequences. For the special case of $K = L$, $f_K = f_L$, and $\alpha = q_L/\pi_L$ this classification is the main result from \cite{DiCapuaKolyvagin}. For $L = \mathbb{Q}_p$ with $\alpha = 1$ and $F_L = x + y + xy$ all such $\alpha\text{-}f_L\text{-norm}$ compatible sequences are interpolated by Coleman series \cite{Coleman1}.
\\
\\
We note that for a series $r(x) \in \mathcal{O}_K[[x]]$ to interpolate some $\alpha\text{-}f_L\text{-norm}$ compatible sequence it is necessary and sufficient for $r$ to satisfy the functional equation
$$r(x) \oplus_L r(x\oplus_K z_1) \oplus_L \ldots \oplus_L r(x\oplus_K z_{q_K-1}) = [\alpha]_L(r([\pi_K]_K(x)))$$

\noindent
The left side of the above equation is guaranteed to converge coefficientwise because $|r(0)|<1$ if $r$ interpolates some $\alpha\text{-}f_L\text{-norm}$ compatible sequence. 
\\
\\
For now we fix some $\alpha \in \mathcal{O}_L$ such that $\alpha \mid q_K/\pi_L$. We define \label{eq:alphamodule} $\mathscr{A}^\alpha_{L,K}$ to be the $\mathcal{O}_L$-module of all power series $r(x) \in \mathcal{O}_K[[x]]$ satisfying $|r(0)| < 1$ and the above functional equation. In the special case $L = K$, $f_K = f_L$, and $\alpha = q_L/\pi_L$ we were able to classify all series in $\mathscr{A}^\alpha_{L,K}$ in \cite{DiCapuaKolyvagin} by mapping $\mathscr{A}^\alpha_{L,K}$ to the kernel of Coleman's trace operator. We use a similar strategy here in order to show that in the more general case $\mathscr{A}^\alpha_{L,K}$ is \say{almost} isomorphic to some submodule of the kernel of Coleman's trace operator. 
\\
\\
It is helpful to define an additional $\mathcal{O}_L$-module first. For $h(x) \in \mathcal{O}_K[[x]]$ consider the functional equation
$$h(x) \oplus_L h(x \oplus_K z_1) \oplus_L \ldots \oplus_L h(x\oplus_K z_{q_K-1}) = 0$$

\noindent
We define $\mathscr{D}_{L,K}$ to be the $\mathcal{O}_L$-module of all series $h(x) \in \mathcal{O}_K[[x]]$ satisfying this functional equation. We define a map $\phi^\alpha_{L,K} : \mathscr{A}^\alpha_{L,K} \rightarrow \mathscr{D}_{L,K}$ as follows: for $r(x) \in \mathscr{A}^\alpha_{L,K}$ we let
$$\phi^\alpha_{L,K}(r) = [q_K/\alpha]_L(r(x))\ominus_L r([\pi_K]_K(x))$$

\noindent
In the above $\ominus_L$ denotes subtraction with respect to the formal group law operation defined by $F_L$. It follows immediately from the definition that $\phi^\alpha_{L,K}$ is a map of $\mathcal{O}_L$-modules. We briefly check that $\phi^\alpha_{L,K}(r) \in \mathscr{D}_{L,K}$ for each series $r \in \mathscr{A}^\alpha_{L,K}$.
\\
\\
In order to check $\phi^\alpha_{L,K}(\mathscr{A}^\alpha_{L,K}) \subseteq \mathscr{D}_{L,K}$ we must consider the following expression and show it is zero:
$$\phi^\alpha_{L,K}(r)(x) \oplus_L \phi^\alpha_{L,K}(r)(x \oplus_K z_1) \oplus_L \ldots \oplus_L \phi^\alpha_{L,K}(r)(x \oplus_K z_{q_K-1})$$

\noindent
For each $i$ we have
$$\phi^\alpha_{L,K}(r)(x \oplus_K z_i) = [q_K/\alpha]_L(r(x \oplus_K z_i)) \ominus r([\pi_K]_K(x))$$

\noindent
Let $\textfrak{F}_0(F_K)$ denote the set of all roots of $f_K(x) = 0$. We let $\sum^{\text{LT},L}$ denote summation with respect to the formal group law operation defined by $F_L$. It then follows from the above that
$$\sum^{\text{LT},L}_{z \in \textfrak{F}_0(F_K)}\phi^\alpha_{L,K}(r)(x\oplus_K z) = \left(\sum^{\text{LT},L}_{z \in \textfrak{F}_0(F_K)}[q_K/\alpha]_L(r(x \oplus_K z))\right) \ominus [q_K]_L(r([\pi_K]_K(x)))$$

\noindent
We also have 
$$\sum^{\text{LT},L}_{z \in \textfrak{F}_0(F_K)} r(x \oplus_K z) = [\alpha]_L(r([\pi_K]_K(x))$$

\noindent
because $r(x) \in \mathscr{A}^\alpha_{L,K}$ which implies that the previous expression must be zero. This completes the proof that $\phi^\alpha_{L,K}(\mathscr{A}^\alpha_{L,K}) \subseteq \mathscr{D}_{L,K}$.
\\
\\
At this point we would like to describe the kernel of the map $\phi^\alpha_{L,K}$. In particular we show that if the kernel is nontrivial then it is generated by a single power series $r_0(x) \in \mathscr{A}^\alpha_{L,K}$ as an $\mathcal{O}_L$-module.
\\
\\
Suppose there exists some nonzero $r(x) \in \mathscr{A}^\alpha_{L,K}$ such that $\phi^\alpha_{L,K}(r) = 0$. This is equivalent to 
$$[q_K/\alpha]_L(r(x)) \ominus_L r([\pi_K]_K(x)) = 0$$

\noindent
Letting $x = 0$ gives 
$$[q_K/\alpha](r(0)) \ominus_L r(0) = 0$$

\noindent
Now $\pi_L \mid q_K/\alpha$ for our choice of $\alpha$, and $|r(0)|<1$. Then the above is only possible if $r(0) = 0$. If $r(u_n) = 0$ for all $n$ then $r(x) = 0$, so we must have some minimal $n = n_0$ such that $r(u_{n_0})$ is nonzero. We get 
$$ [q_K/\alpha](r(u_{n_0})) = r(u_{n_0-1}) = 0$$

\noindent
where we take $u_{n_0-1} = 0$ if $n_0 = 0$. This implies $r(u_{n_0})$ is a torsion point of $F_L$. At this point we let $\textfrak{F}_n(F_L)$ denote the set of roots of $f_L^{(n+1)}(x) = 0$. Let $\textfrak{F}_\infty(F_L)$ denote the set of all torsion points of $F_L$. The above implies the sequence $(r(u_n))_{n \ge 0}$ is such that $r(u_n) \in \textfrak{F}_\infty(F_L)$ for each $n \ge 0$ and such that $[q_K/\alpha](r(u_{n+1})) = r(u_n)$ for each $n \ge 0$.
\\
\\
For each nonzero series $r(x) \in \mathscr{A}^\alpha_{L,K}$ satisfying $\phi^\alpha_{L,K}(r) = 0$ we define the first index $n_r$ of $r$ to be the smallest integer $n$ such that $r(u_n)$ is nonzero. We define the second index $m_r$ of $r$ to be the minimal $m$ such that $r(u_{n_r}) \in \textfrak{F}_m(F_L)$. Note that there is an upper bound on all possible values of $m_r$ because $r(u_{n_r})$ must be a root of $[q_K/\alpha](x) = 0$.
\\
\\
Out of all possible series $r(x)$ in the kernel of $\phi^\alpha_{L,K}$ consider the subset of $r(x)$ such that $n_r$ is minimal. Let $r_0(x)$ be any series in this subset with maximal second index $m_r$. We show that this choice of $r_0(x)$ implies that $r_0(x)$ generates the kernel of $\phi^\alpha_{L,K}$ as an $\mathcal{O}_L$-module.
\\
\\
Let $r$ be any series $r(x) \in \mathscr{A}^\alpha_{L,K}$ satisfying $\phi^\alpha_{L,K}(r) = 0$. We show there exists a unique $a_r \in \mathcal{O}_L$ such that $[a_r]_L(r_0(x)) = r(x)$. For any index $N > n_{r_0}$ we have $r_0(u_N)$ and $r(u_N)$ are both torsion points of $F_L$. Let $m_N$ be the smallest index such that $\textfrak{F}_{m_N}(F_L)$ contains $r_0(u_N)$. Then $r_0(u_N)$ generates $\textfrak{F}_{m_N}(F_L)$ as an $\mathcal{O}_L$-module. We must also have $r(u_N) \in \textfrak{F}_{m_N}(F_L)$ because of our choice of $r_0(x)$. Therefore there exists some $a_N \in \mathcal{O}_L$ such that $[a_N](r_0(u_N)) = r(u_N)$. This implies that $[a_N](r_0(u_n)) = r(u_n)$ for each $n \le N$ because $(r_0(u_n))$ and $(r(u_n))$ are both compatible with respect to $[q_K/\alpha](x)$. 
\\
\\
The above equalities imply $\lim_{N\rightarrow \infty}a_N$ exists in $\mathcal{O}_L$ and if we take $a_r$ to be this limit we will have $[a_r](r_0(u_n)) = r(u_n)$ for all $n \ge 0$. Equality on all torsion points of $F_K$ implies $[a_r](r_0(x)) = r(x)$ as series. For a proof of this see the "Uniqueness Principle" of \cite{Coleman1}. This completes the proof that if $\phi^\alpha_{L,K}$ has nontrivial kernel on $\mathscr{A}^\alpha_{L,K}$ then the kernel is generated by a single series.
\\
\\
In order to show that $\mathscr{A}^\alpha_{L,K}$ is \say{almost} isomorphic to a submodule of the kernel of Coleman's trace operator we define some additional $\mathcal{O}_L$-modules. Let $\mathscr{D}_{L,K,\pi_L}$ be the submodule of $\mathscr{D}_{L,K}$ consisting of all series $h(x) \in \mathscr{D}_{L,K}$ such that $\pi_L \mid h(x)$ in $\mathcal{O}_K[[x]]$. Let $\mathscr{L}_{F_K}$ denote Coleman's trace operator with respect to $F_K$. Recall that $\mathscr{L}_{F_K}$ is the operator defined on $\mathcal{O}_K[[x]]$ by $\mathscr{L}_{F_K}(f)([\pi_K]_K(x)) = \sum_{z\in \textfrak{F}_0(F_K)} f(x \oplus_K z)$.
\\
\\
We define $\mathscr{C}_{L,K}$ to be the $\mathcal{O}_L$-module of power series $g(x) \in \pi_L\mathcal{O}_K[[x]]$ such that $\mathscr{L}_{F_K}(g) = 0$. We let $\log_{F_L}$ and $\exp_{F_L}$ denote the logarithm and the exponential of $F_L$ respectively. Then one can show that $\log_{F_L} : \mathscr{D}_{L,K,\pi_L} \rightarrow \mathscr{C}_{L,K}$ is an isomorphism of $\mathcal{O}_L$-modules with inverse $\exp_{F_L}$.
\\
\\
\hypertarget{bottompg5}{We follow \cite{DiCapuaKolyvagin} to show that $\log_{F_L} : \mathscr{D}_{L,K,\pi_L} \rightarrow \mathscr{C}_{L,K}$ is an isomorphism of $\mathcal{O}_L$-modules. It is well known that $\log_{F_L} : \pi_L\mathcal{O}_K \rightarrow \pi_L\mathcal{O}_K$ and $\exp_{F_L} : \pi_L\mathcal{O}_K \rightarrow \pi_L\mathcal{O}_K$ are inverse isomorphisms of $\mathcal{O}_L$-modules. See for example Proposition 7.17 and Proposition 2.4 in \cite{Kolyvagin}, agreeing with the $\mathcal{O}_L$-action follows if we consider Theorem 2 in section 5.1 of \cite{BorevichShafarevich}. The same estimates of divisibility of $\log_{F_L}(a)$ and $\exp_{F_L}(b)$ depending on divisibility of $a, b$ in the above proof imply that $\log_{F_L} : \pi_L\mathcal{O}_K[[x]] \rightarrow \pi_L\mathcal{O}_K[[x]]$ and $\exp_{F_L} : \pi_L\mathcal{O}_K[[x]] \rightarrow \pi_L\mathcal{O}_K[[x]]$ are defined as coefficientwise limits. The remaining claims follow because the series are free for substitutions $x \in \pi_L\mathcal{O}_K$ and coefficientwise limits agree with composition of functions on $\pi_L\mathcal{O}_K$.}
\\
\\
Recall we have a map of $\mathcal{O}_L$-modules $\phi^\alpha_{L,K} : \mathscr{A}^\alpha_{L,K} \rightarrow \mathscr{D}_{L,K}$. We now need the following lemma:

\begin{lemma}
    \label{lem:2.0.1}
    If $h(x) \in \mathscr{D}_{L,K}$ and $\pi_K \mid h(0)$ then $\pi_K \mid h(x)$ in $\mathcal{O}_K[[x]]$.
\end{lemma}

\noindent
Proof: suppose there exists some $h(x) \in \mathscr{D}_{L,K}$ such that $\pi_K \mid h(0)$ and such that $h(x) \not\equiv 0 \mod \pi_K$. We know $h$ satisfies
$$\sum^{\text{LT},L}_{z\in\textfrak{F}_0(F_K)}h(x\oplus_K z) = 0$$

\noindent
Considering the above equation mod $z_1$ gives
$$[q_K](h(x)) \equiv 0 \mod z_1$$

\noindent
Then there is some unit $u \in \mathcal{O}_L$ and some integer $n$ such that $q_K = u\pi_L^n$ in $\mathcal{O}_L$. It follows that
$$[q_K](h(x)) \equiv [u\pi_L^n](h(x)) \equiv [u](h(x)^{q_L^n}) \mod z_1$$

\noindent
Now $[u](h(x)^{q_L^n}) \not\equiv 0 \mod z_1$ because $h(0) \equiv 0 \mod z_1$ and because $h(x)^{q_L^n} \not\equiv 0 \mod \pi_K$. It follows that $[u](h(x)^{q_L^n})$ cannot be zero mod $z_1$ which contradicts 
$$[q_K](h(x)) \equiv 0 \mod z_1$$

\noindent
We conclude if $h(x) \in \mathscr{D}_{L,K}$ and $\pi_K \mid h(0)$ then $\pi_K \mid h(x) \in \mathcal{O}_K[[x]]$.
\\
\\
If $r(x) \in \mathscr{A}^\alpha_{L,K}$ then 
$$\phi^\alpha_{L,K}(r) = [q_K/\alpha]_L(r(x)) \ominus_L r([\pi_K]_K(x))$$

\noindent
It follows that 
$$\phi^\alpha_{L,K}(r)(0) = [q_K/\alpha - 1]_L(r(0))$$

\noindent
From the above we conclude that for every $h(x) \in \phi^\alpha_{L,K}(\mathscr{A}^\alpha_{L,K})$ we have $\pi_K \mid h(0)$. Then we also have $\pi_K \mid h(x)$ in $\mathcal{O}_K[[x]]$ by \hyperref[lem:2.0.1]{lemma 2.1}.
\\
\\
At this point we fix some integer $r$ sufficiently large so that $[\pi_L^r](\pi_K\mathcal{O}_K[[x]]) \subseteq \pi_L\mathcal{O}_K[[x]]$. Note that $[\pi_L^r] : \pi_K\mathcal{O}_K[[x]] \rightarrow \pi_L\mathcal{O}_K[[x]]$ defined by sending $f(x)$ to $[\pi_L^r](f(x))$ is an injective map of $\mathcal{O}_L$-modules since $[\pi_L](x)$ has an inverse under composition of formal power series in $L[[x]]$. 
\\
\\
It follows from the above that the composition $[\pi_L^r] \circ \phi^\alpha_{L,K} : \mathscr{A}^\alpha_{L,K} \rightarrow \mathscr{D}_{L,K,\pi_L}$ is a map of $\mathcal{O}_L$-modules with kernel equal to the kernel of $\phi^\alpha_{L,K}$. We then compose this map with the logarithm of $F_L$ to get the following:

\begin{theorem}
    \label{thm:2.0.2}
    The map $\mathscr{A}^\alpha_{L,K} \rightarrow \mathscr{C}_{L,K}$ defined by sending $r(x) \in \mathscr{A}^\alpha_{L,K}$ to $\log_{F_L}([\pi_L^r](\phi^\alpha_{L,K}(r))) \in \mathscr{C}_{L,K}$ either has trivial kernel or has a kernel generated by a single series as an $\mathcal{O}_L$-module. In this sense the $\mathcal{O}_L$-module $\mathscr{A}^\alpha_{L,K}$ is \say{almost} isomorphic to some submodule of the kernel of Coleman's trace operator.
\end{theorem}

\noindent
We will soon be able to show that for many values of $\alpha$ the kernel mentioned in the above theorem must be trivial. We begin by asking the following question: for which values of $\alpha \in \mathcal{O}_L$ do there exist series $r(x) \in \mathscr{A}^\alpha_{L,K}$ such that $r(x)$ does not vanish mod $\pi_K$?
\\
\\
We will show that for a fixed choice of $L$ and $K$ there is exactly one choice of $\alpha$ such that the module $\mathscr{A}_{L,K}^\alpha$ can contain series $r(x)$ such that $r(x) \not\equiv 0 \mod \pi_K$. For the remaining choices of $\alpha$, we will show that all series in $\mathscr{A}^\alpha_{L,K}$ are coming from the $\alpha$-eigenspace of Coleman's trace operator on $\pi_L\mathcal{O}_K[[x]]$ in some precise sense.
\\
\\
For a fixed choice of $L$ and $K$ we first suppose $\mathscr{A}^\alpha_{L,K}$ contains some series $r(x)$ such that $r(x)$ does not vanish mod $\pi_K$. Recall this means $r(x)$ is a solution to the functional equation
$$r(x) \oplus_L r(x\oplus_K z_1) \oplus_L \ldots \oplus_L r(x\oplus_K z_{q_K-1}) = [\alpha]_L(r([\pi_K]_K(x)))$$

\noindent
Considering the above equation mod $z_1$ we get that
$$[q_K](r(x)) \equiv [\alpha]_L(r(x^{q_K})) \mod z_1$$

\noindent
Since $K$ is a finite extension of $L$ we have there exists a positive integer $f_{K/L}$ such that $q_K = q_L^{f_{K/L}}$. This implies
$$[q_K](r(x)) \equiv [\alpha]_L(r(x^{q_K})) \equiv [\alpha \pi_L^{f_{K/L}}](r(x)) \mod z_1$$

\noindent
The above is equivalent to
$$[q_K - \alpha\pi_L^{f_{K/L}}](r(x)) \equiv 0 \mod z_1$$

\noindent
However, note that if $q_K - \alpha\pi_L^{f_{K/L}}$ is any nonzero element of $\mathcal{O}_L$, then $[q_K - \alpha\pi_L^{f_{K/L}}](r(x))$ must be nonzero mod $z_1$ because $\pi_K \mid r(0)$ and because $r(x) \not\equiv 0 \mod z_1$. This implies that if $\alpha \in \mathcal{O}_L$ is such that $\mathscr{A}^\alpha_{L,K}$ contains series which do not vanish mod $\pi_K$, then we must have $\alpha = q_K/\pi_L^{f_{K/L}}$. Note for all other values of $\alpha$, we must have that the kernel of the map from \hyperref[thm:2.0.2]{Theorem 2.2} is trivial (to see this note that if $r \neq 0$ is in the kernel of $\phi^\alpha_{L,K}$ then the sequence $(r(u_n))$ is a sequence of torsion points. The levels of the points in this sequence is not bounded).
\\
\\
For all values of $\alpha$ such that $\alpha \neq q_K/\pi_L^{f_{K/L}}$ we show that all series in $\mathscr{A}^\alpha_{L,K}$ are coming from the $\alpha$-eigenspace of $\mathscr{L}_{F_K}$ in $\pi_L\mathcal{O}_K[[x]]$ in a certain sense. In order to prove this we define the $\mathcal{O}_L$-module $\mathscr{A}^\alpha_{L,K,\pi_L}$ to be the submodule of $\mathscr{A}^\alpha_{L,K}$ consisting of all series $r(x) \in \mathscr{A}^\alpha_{L,K}$ such that $\pi_L \mid r(x)$ in $\mathcal{O}_K[[x]]$. For all values of $\alpha \neq q_K/\pi_L^{f_{K/L}}$ we have that $\pi_K \mid r(x)$ for each $r(x) \in \mathscr{A}^\alpha_{L,K}$. It follows that we can fix a sufficiently large integer $s$ such that $[\pi_L^s](\mathscr{A}^\alpha_{L,K}) \subseteq \mathscr{A}^\alpha_{L,K,\pi_L}$.
\\
\\
Next we show $\mathscr{A}^\alpha_{L,K,\pi_L}$ is isomorphic to the $\alpha$-eigenspace of $\mathscr{L}_{F_K}$ in $\pi_L\mathcal{O}_K[[x]]$. We define $\mathscr{E}^\alpha_{L,K}$ to be the $\mathcal{O}_L$-module of all series $g(x) \in \pi_L\mathcal{O}_K[[x]]$ such that
$$\mathscr{L}_{F_K}(g) = \alpha g$$

\noindent
Note that writing the above equality is equivalent to writing $g$ satisfies
$$\sum_{z\in\textfrak{F}_0(F_K)}g(x\oplus_K z) = \alpha g([\pi_K]_K(x))$$

\noindent
We check that the map $\log_{F_L} : \mathscr{A}^\alpha_{L,K,\pi_L} \rightarrow \mathscr{E}^\alpha_{L,K}$ defined by sending $h(x) \in \mathscr{A}^\alpha_{L,K,\pi_L}$ to the composition $\log_{F_L}(h(x))$ is an isomorphism of $\mathcal{O}_L$-modules with inverse $\exp_{F_L}$.
\\
\\
First it is useful to note that if $h(x) \in \mathscr{A}^\alpha_{L,K,\pi_L}$ then $\pi_L \mid h(x)$ so that the composition $\log_{F_L}(h(x))$ will also live in $\pi_L\mathcal{O}_K[[x]]$ by the arguments on \hyperlink{bottompg5}{page 10} of this paper. We must check that if $h(x) \in \mathscr{A}^\alpha_{L,K,\pi_L}$ then the composition $\log_{F_L}(h(x))$ satisfies the correct functional equation. 
\\
\\
For $h(x)$ to live in $\mathscr{A}^\alpha_{L,K,\pi_L}$ means $h(x)$ satisfies
$$\sum^{\text{LT},L}_{z \in \textfrak{F}_0(F_K)} h(x \oplus_K z) = [\alpha]_L(h([\pi_K]_K(x))$$

\noindent
Applying $\log_{F_L}$ to both sides of this equation implies
$$\sum_{z \in \textfrak{F}_0(F_K)} \log_{F_L}(h(x \oplus_K z)) = \alpha\log_{F_L}(h([\pi_K]_K(x)))$$

\noindent
We get the left side of the above equation because $\log_{F_L}$ takes summation with respect to the formal group law $F_L$ to an ordinary summation. We get the right side of the above equation because we have the identity $\log_{F_L}([\alpha]_L(x)) = \alpha \log_{F_L}(x)$ for each $\alpha \in \mathcal{O}_L$. Because the above equation is exactly the equation defining $\mathscr{E}^\alpha_{L,K}$ we conclude $\log_{F_L}(\mathscr{A}^\alpha_{L,K,\pi_L}) \subseteq \mathscr{E}^\alpha_{L,K}$.
\\
\\
We also need to check that $\exp_{F_L}(\mathscr{E}^\alpha_{L,K}) \subseteq \mathscr{A}^\alpha_{L,K,\pi_L}$. Once we have this inclusion it will follow that $\log_{F_L} : \mathscr{A}^\alpha_{L,K,\pi_L} \rightarrow \mathscr{E}^\alpha_{L,K}$ is an isomorphism of $\mathcal{O}_L$-modules with inverse $\exp_{F_L}$ because $\log_{F_L}(x)$ and $\exp_{F_L}(x)$ are inverses under composition.
\\
\\
Recall $g(x) \in \mathscr{E}^\alpha_{L,K}$ exactly when $\pi_L \mid g(x) \in \mathcal{O}_K[[x]]$ and $g$ satisfies the functional equation
$$\sum_{z\in\textfrak{F}_0(F_K)}g(x\oplus_K z) = \alpha g([\pi_K]_K(x))$$

\noindent
Again by the arguments on \hyperlink{bottompg5}{page 10} of this paper we have $g(x) \in \pi_L\mathcal{O}_K[[x]]$ implies the composition $\exp_{F_L}(g(x))$ lives in $\pi_L\mathcal{O}_K[[x]]$. Applying $\exp_{F_L}$ to both sides of the above equation gives
$$\sum^{\text{LT},L}_{z\in\textfrak{F}_0(F_K)}\exp_{F_L}(g(x\oplus_K z)) = [\alpha]_L (\exp_{F_L}(g([\pi_K]_K(x))))$$

\noindent
We get the left side of the above equation because $\exp_{F_L}$ takes summation over ordinary addition to summation over addition with respect to the formal group law $F_L$. We get the right side of the above equation because we have the identity $\exp_{F_L}(\alpha x) = [\alpha]_L(\exp_{F_L}(x))$ for each $\alpha \in \mathcal{O}_L$. Because the above equation is exactly the equation defining $\mathscr{A}^\alpha_{L,K}$ we conclude $\exp_{F_L}(\mathscr{E}^\alpha_{L,K}) \subseteq \mathscr{A}^\alpha_{L,K,\pi_L}$. This completes the proof of the following:

\begin{theorem}
    \label{thm:2.0.3}
    $\log_{F_L}: \mathscr{A}^\alpha_{L,K,\pi_L} \rightarrow \mathscr{E}^\alpha_{L,K}$ is an isomorphism of $\mathcal{O}_L$-modules with inverse $\exp_{F_L}$.
\end{theorem}

\noindent
We would now like to describe $\mathscr{E}^\alpha_{L,K}$ for each $\alpha \in \mathcal{O}_L$. This was already essentially done in \cite{DiCapuaKolyvagin}. For $\alpha$ a unit in $\mathcal{O}_L$ we have that $\mathscr{E}^\alpha_{L,K}$ is empty. This follows from lemma 6 of \cite{Coleman1}. For the remaining cases where $\pi_L \mid \alpha$ in $\mathcal{O}_L$ we closely follow the arguments starting on page 44 of \cite{DiCapuaKolyvagin}.
\\
\\
For $\pi_L \mid \alpha$ our goal is to show $\mathscr{E}^\alpha_{L,K}$ is isomorphic to to $\mathscr{C}_{L,K}$ as $\mathcal{O}_L$-modules. Recall that $\mathscr{C}_{L,K}$ is defined to be the $\mathcal{O}_L$-module of all power series $h(x) \in \pi_L\mathcal{O}_K[[x]]$ such that $\mathscr{L}_{F_K}(h) = 0$. We need to define two series. $k(x)$ is defined to be any series $k(x) \in x^{q-1}\mathcal{O}_K[[x]]$ such that $k(u_0) = \pi_K/(q-1)$. For the existence of such $k(x)$ see page 19 of \cite{DiCapuaKolyvagin}. We also define $w(x)$ to be the unique series $w(x) \in \mathcal{O}_K[[x]]$ such that $\mathscr{L}_{F_K}(k) = \pi_Kw$. The series $k$ and $w$ originally appear in \cite{Coleman2}. 
\\
\\
We define a map of $\mathcal{O}_L$-modules $\rho_\alpha : \mathscr{E}^\alpha_{L,K} \rightarrow \mathscr{C}_{L,K}$ as follows:
$$\rho_\alpha(g) = g(x) - \frac{\alpha k(x) g([\pi_K]_K(x))}{\pi_Kw([\pi_K]_K(x))}$$

\noindent
We are able to prove the following:
\begin{theorem}
    $\rho_\alpha : \mathscr{E}^\alpha_{L,K} \rightarrow \mathscr{C}_{L,K}$ is an isomorphism of $\mathcal{O}_L$-modules.
\end{theorem}

\noindent
We must first check that $\rho_\alpha(g)$ is contained in $\mathscr{C}_{L,K}$ for each $g(x) \in \mathscr{E}^\alpha_{L,K}$.
\\
\\
Note that if $\pi_L \mid g(x)$ in $\mathcal{O}_K[[x]]$ then we also have $\pi_L \mid \rho_\alpha(g)$ in $\mathcal{O}_K[[x]]$. This is because $\pi_L \mid \alpha$, $\pi_L \mid g([\pi_K]_K(x))$, and $w([\pi_K]_K(x))$ is a unit in $\mathcal{O}_K[[x]]$. Then once we check that $\mathscr{L}_{F_K}(\rho_\alpha(g)) = 0$ we can conclude $\rho_\alpha(g) \in \mathscr{C}_{L,K}$.
\\
\\
Consider the expression 
$$\mathscr{L}_{F_K}(\rho_\alpha(g)) = \mathscr{L}_{F_K}(g) - \mathscr{L}_{F_K}(\frac{\alpha k(x) g([\pi_K]_K(x))}{\pi_Kw([\pi_K]_K(x))})$$

\noindent
We have $\mathscr{L}_{F_K}(g) = \alpha g$ because $g(x) \in \mathscr{E}^\alpha_{L,K}$. By the construction of the series $k$ and $w$ and the arguments in \cite{Coleman2} we know that
$$\mathscr{L}_{F_K}(\frac{k(x)g([\pi_K]_K(x))}{\pi_Kw([\pi_K]_K(x))}) = g(x)$$

\noindent
Combining this with the previous expression for $\mathscr{L}_{F_K}(\rho_\alpha(g))$ gives us
$$\mathscr{L}_{F_K}(\rho_\alpha(g)) = \alpha g - \alpha g = 0$$

\noindent
so we conclude $\rho_\alpha(g) \in \mathscr{C}_{L,K}$. Next we would like to show the map $\rho_\alpha$ is injective. In section 3.1 of \cite{DiCapuaKolyvagin} we showed the analog of this map is injective by assuming there is a maximal $N$ such that $x^N \mid g(x)$ and arriving at a contradiction. The same argument works here to show that $\rho_\alpha$ is injective. 
\\
\\
We note that if we further restrict $\alpha \in \mathcal{O}_L$ by stipulating $\pi_L^2 \mid \alpha$ or if we assume the extension $K/L$ is ramified, so that $|\pi_L/\pi_K| < 1$, then we have an additional proof that $\rho_\alpha$ is injective. For this proof let $M$ be the largest integer such that $\pi_K^M \mid g(x)$ in $\mathcal{O}_K[[x]]$. $M$ must exist if $g(x) \neq 0$. We are assuming $\rho_\alpha(g) = 0$, so we have
$$\rho_\alpha(g) = g(x) - \frac{\alpha k(x) g([\pi_K]_K(x))}{\pi_Kw([\pi_K]_K(x))} = 0$$

\noindent
We note that $\pi_K^M \mid g([\pi_K]_K(x))$ and $|\alpha/\pi_K| < 1$, so we must have $\pi_K \mid \alpha/\pi_K$. This implies that $\pi_K^{M+1} \mid g(x)$ which is a contradiction. We conclude the map $\rho_\alpha$ is injective. 
\\
\\
Finally we must show $\rho_\alpha : \mathscr{E}^\alpha_{L,K} \rightarrow \mathscr{C}_{L,K}$ is surjective, and then we can conclude $\rho_\alpha$ is an isomorphism of these $\mathcal{O}_L$-modules. Note for all $\alpha \in \mathcal{O}_L$ such that $\pi_L \mid \alpha$, the proof used in section 3.1 of \cite{DiCapuaKolyvagin} used to show $\rho_\lambda : N_\lambda \rightarrow M$ is surjective also works to show $\rho_\alpha : \mathscr{E}^\alpha_{L,K} \rightarrow \mathscr{C}_{L,K}$ is surjective. Here we impose the additional hypothesis $|\alpha/\pi_K| < 1$ in order to give a different proof that $\rho_\alpha : \mathscr{E}^\alpha_{L,K} \rightarrow \mathscr{C}_{L,K}$ is surjective in this special case.
\\
\\
Suppose $|\alpha/\pi_K| < 1$. We show that $\rho_\alpha : \mathscr{E}^\alpha_{L,K} \rightarrow \mathscr{C}_{L,K}$ is surjective. We take arbitrary $h(x) \in \mathscr{C}_{L,K}$. We construct a series $g(x)$ such that $\rho_\alpha(g) = h$.
\\
\\
We define a sequence of series $g_i(x) \in \pi_L\mathcal{O}_K[[x]]$ such that the infinite sum $\sum_{i=1}^\infty g_i(x)$ converges coefficientwise to some series $g(x) \in \pi_L\mathcal{O}_K[[x]]$. We will construct the sequence $(g_i)$ in such a way that the partial sums $\rho_\alpha( \sum_{i=1}^N g_i(x))$ converge coefficientwise to $h(x)$. Take $g_1(x) = h(x)$, so we have
$$\rho_\alpha(g_1(x)) = h(x) - \frac{\alpha k(x) h([\pi_K]_K(x))}{\pi_Kw([\pi_K]_K(x))}$$

\noindent
Take $g_2(x) = h(x) - \rho_\alpha(g_1(x))$, so that
$$g_2(x) = \frac{\alpha k(x) h([\pi_K]_K(x))}{\pi_K w([\pi_K]_K(x))}$$

\noindent
If $M_h$ is the largest integer such that $\pi_K^{M_h} \mid h(x)$ we note that $\pi_K^{M_h + 1} \mid g_2(x)$. We also note that
$$\rho_\alpha(g_1 + g_2) = \rho_\alpha(g_1) + \rho_\alpha(g_2) = h(x) - \frac{\alpha k(x) g_2([\pi_K]_K(x))}{\pi_Kw([\pi_K]_K(x))}$$

\noindent
implying that $\pi_K^{M_h + 2} \mid h(x) - \rho_\alpha(g_1 + g_2)$. 
\\
\\
Take some integer $N \ge 2$ and suppose for all $1 \le n \le N$ we have found series $g_n$ such that $\pi_K^{M_h+n-1} \mid g_n(x)$ and such that 
$$\pi_K^{M_h+n} \mid h(x) - \rho_\alpha(\sum_{i=1}^ng_i(x))$$

\noindent
We show it is possible to pick $g_{N+1}(x)$ such that $\pi_K^{M_h + N} \mid g_{N+1}(x)$ and such that 
$$\pi_K^{M_h + N +1} \mid h(x) - \rho_\alpha(\sum_{i=1}^{N+1}g_i(x))$$

\noindent
Take $g_{N+1}(x) = h(x) - \rho_\alpha(\sum_{i=1}^Ng_i(x))$. For this choice of $g_{N+1}(x)$ we immediately get $\pi_K^{M_h + N} \mid g_{N+1}(x)$. We also have
$$\rho_\alpha(g_{N+1}(x)) = h(x) - \rho_\alpha(\sum_{i=1}^N g_i(x)) - \frac{\alpha k(x) g_{N+1}([\pi_K]_K(x))}{\pi_Kw([\pi_K]_K(x))}$$

\noindent
implying $\pi_K^{M_h + N + 1} \mid h(x) - \rho_\alpha(\sum_{i=1}^{N+1}g_i(x))$.
\\
\\
For the sequence $(g_i(x))$ constructed above we have $\pi_K^{M_h + i -1} \mid g_i(x)$ for each $i$ implying that the partial sums $\sum_{i=1}^N g_i(x)$ converge coefficientwise to some $g(x) \in \pi_L\mathcal{O}_K[[x]]$. It follows that the expressions $\rho_\alpha(\sum^N_{i=1}g_i(x))$ also converge coefficientwise to $\rho_\alpha(g)$. We have already seen these expressions converge coefficientwise to $h(x)$, so we conclude that $\rho_\alpha(g) = h$.
\\
\\
This completes the proof that $\rho_\alpha : \mathscr{E}^\alpha_{L,K} \rightarrow \mathscr{C}_{L,K}$ is an isomorphism of $\mathcal{O}_L$-modules. One can use the methods from section 2.2 of \cite{DiCapuaKolyvagin} to describe all series in $\mathscr{C}_{L,K}$. Recall $\log_{F_L} : \mathscr{A}^\alpha_{L,K,\pi_L} \rightarrow \mathscr{E}^\alpha_{L,K}$ is also an isomorphism of $\mathcal{O}_L$-modules. We also mentioned earlier that for all values of $\alpha \neq q_K/\pi_L^{f_{K/L}}$ it is possible to fix some integer $s$ such that $[\pi^s](\mathscr{A}^\alpha_{L,K}) \subseteq \mathscr{A}^\alpha_{L,K,\pi_L}$. This completes our description of $\mathscr{A}^\alpha_{L,K}$ for all values of $\alpha \neq q_K/\pi_L^{f_{K/L}}$ which shows that in a certain sense all series in $\mathscr{A}^\alpha_{L,K}$ are coming from the $\alpha$-eigenspace of $\mathscr{L}_{F_K}$ in $\pi_L\mathcal{O}_K[[x]]$.
\\
\\
At this point we can also say more about the kernel of the map $\phi^\alpha_{L,K} : \mathscr{A}^\alpha_{L,K} \rightarrow \mathscr{D}_{L,K}$ from \hyperref[thm:2.0.2]{Theorem 2.2} for all values of $\alpha \neq q_K/\pi_L^{f_{K/L}}$ and such that $\alpha \mid q_K/\pi_L$. Suppose $r(x)$ is some nonzero series in $\mathscr{A}^\alpha_{L,K}$ such that $\phi^\alpha_{L,K}(r) = 0$. This is equivalent to writing
$$[q_K/\alpha]_L(r(x)) \ominus_L r([\pi_K]_K(x)) = 0$$

\noindent 
Letting $x = 0$ gives $[q_K/\alpha - 1](r(0)) = 0$ which is only possible if $r(0) = 0$ because $q_K/\alpha - 1$ is a unit and $\pi_K \mid r(0)$. The above implies the sequence $(r(u_n))$ satisfies each $r(u_n) \in \textfrak{F}_\infty(F_L)$ and 
$$[q_K/\alpha]_L(r(u_{n+1})) = r(u_n)$$

\noindent
for each $n \ge 0$. Because $\pi_L \mid q_K/\alpha$ the above can only happen if $\lim_{n\rightarrow \infty}|r(u_n)| = 1$. However we know $\pi_K \mid r(x)$ for every $r(x) \in \mathscr{A}^\alpha_{L,K}$ whenever $\alpha$ is not $q_K/\pi_L^{f_{K/L}}$. $\pi_K \mid r(x)$ implies $|r(u_n)| \le |\pi_K|$ for each $n$, so we get a contradiction if the kernel of $\phi^\alpha_{L,K}$ is nontrivial in this case. For such values of $\alpha$ we conclude that the map taking $r(x) \in \mathscr{A}^\alpha_{L,K}$ to the series $\log_{F_L}([\pi_L^r](\phi^\alpha_{L,K}(r))) \in \mathscr{C}_{L,K}$ is an injective map of $\mathcal{O}_L$-modules. In other words, the map from \hyperref[thm:2.0.2]{Theorem 2.2} is an isomorphism between $\mathscr{A}^\alpha_{L,K}$ and some submodule of $\mathscr{C}_{L,K}$ in these cases.

\subsection{The special case $\alpha = q_K/\pi_L^{f_{K/L}}$}

For now we fix $\alpha = q_K/\pi_L^{f_{K/L}}$. As we have already seen this is the only value of $\alpha$ for which there can exist $r(x) \in \mathscr{A}^\alpha_{L,K}$ such that $r(x)$ does not vanish mod $\pi_K$. Our goal in this section is to give a criterion for checking the following condition: given a series $s(x) \in \mathcal{O}_K[[x]]$ whose coefficients are all units does there exist some $r(x) \in \mathscr{A}^\alpha_{L,K}$ such that $r(x) \equiv s(x) \mod \pi_K$?
\\
\\
For now we fix $f_L(x)$ to be the series $f_L(x) = \pi_Lx + x^{q_L}$. We also fix $f_K(x)$ to be $f_K(x) = \pi_Kx + x^{q_K}$. We also assume $q_L \ge 3$ and $\pi_K^2 \mid \pi_L$. We are able to prove the following under these conditions:

\begin{theorem}
\label{thm:2.1.1}
Assume $\pi_L \mid \alpha = q_K/\pi_L^{f_{K/L}}$. Then for any series $s(x) \in \mathcal{O}_K[[x]]$ whose coefficients are all units, we have there exists $r(x) \in \mathscr{A}^\alpha_{L,K}$ such that $r(x) \equiv s(x) \mod \pi_K$ if and only if
$$\sum^{\text{LT},L}_{z \in \textfrak{F}_0(F_K)} s(x \oplus_K z) \equiv [\alpha]_L(s([\pi_K]_K(x)) \mod \pi_K^2$$

\noindent

\end{theorem}

\noindent
We will need the following lemma:

\begin{lemma}
    \label{lem:2.1.2}
    The formal group law $F_L(x,y)$ associated to $f_L(x) = \pi_Lx + x^{q_L}$ satisfies
    \[
    F_L(x,y) \equiv x + y \mod \deg q_L
    \]
\end{lemma}
\noindent
where equivalence mod deg $q_L$ means that $F_L(x,y) - (x+y)$ lives in the ideal of $\mathcal{O}_K[[x,y]]$ generated by the series $x^{q_L}, x^{q_L-1}y, x^{q_L-2}y^2, \ldots y^{q_L}$.
\\
\\
Proof: by section 3.5 of chapter 6 of \cite{CasselsFrohlich} we have that $F_L$ is the unique power series in $\mathcal{O}_K[[x,y]]$ satisfying 
$$F_L(x,y) \equiv x + y \mod \deg 2$$

\noindent
and
$$f_L(F_L(x,y)) = F_L(f_L(x),f_L(y))$$

\noindent
The existence of $F_L$ is proved by Proposition 5 in section 3.5 of chapter 6 of \cite{CasselsFrohlich}. Proposition 5 also describes how to recursively solve for the coefficients of $F_L(x,y)$ using $f_L(x)$. We use this method to prove \hyperref[lem:2.1.2]{Lemma 2.6}.
\\
\\
Using the language of Propsition 5 from \cite{CasselsFrohlich} we take $\phi^{(1)} = \phi_1 = x + y$. We have the equation
$$f_L\circ \phi^{(1)} \equiv \phi^{(1)}\circ f_L + E_2 \mod \deg 3$$

\noindent
If we show that the difference $f_L\circ\phi^{(1)} - \phi^{(1)}\circ f_L$ is actually zero mod degree $q_L$ this proves \hyperref[lem:2.1.2]{Lemma 2.6}. First we compute $f_L\circ \phi^{(1)} = f_L(x+y)$.
$$f_L(x + y) = \pi_L(x+y) + (x+y)^{q_L} = \pi_Lx + \pi_Ly + x^{q_L} + y^{q_L} + \sum_{k=1}^{q_L - 1} {q_L \choose k} x^ky^{q_L - k}$$

\noindent
Next we have
$$\phi^{(1)}\circ f_L = f_L(x) + f_L(y) = \pi_Lx + x^{q_L} + \pi_Ly + y^{q_L}$$

\noindent
The above computations imply the difference $f_L\circ \phi^{(1)} - \phi^{(1)}\circ f_L$ is congruent to zero mod degree $q_L$. This implies for each index $i$ with $2 \le i < q_L$ we have $E_i$ as defined in the proof of Proposition 5 will be zero. Since $E_i = 0$ implies $\phi_i = 0$ we conclude from the proof of Proposition 5 that $F_L(x,y) \equiv x + y \mod \deg q_L$.
\\
\\
We now move to the proof of \hyperref[thm:2.1.1]{Theorem 2.5}. First suppose for some $s(x) \in \mathcal{O}_K[[x]]$ such that all of the coefficients of $s$ are units we have there exists some $r(x) \in \mathscr{A}^\alpha_{L,K}$ such that $r(x) \equiv s(x) \mod \pi_K$. We show this implies
$$\sum^{\text{LT},L}_{z \in \textfrak{F}_0(F_K)} s(x \oplus_K z) \equiv [\alpha]_L(s([\pi_K]_K(x)) \mod \pi_K^2$$

\noindent
Define the operator $T$ on power series in $\mathcal{O}_K[[x]]$ which have constant term divisible by $\pi_K$ by sending $f \in \mathcal{O}_K[[x]]$ (with $|f(0)| < 1$) to the series
$$T(f) = \left(\sum^{\text{LT},L}_{z \in \textfrak{F}_0(F_K)} f(x \oplus_K z) \right) \ominus_L [\alpha]_L(f([\pi_K]_K(x))$$

\noindent
Note that $f \in \mathscr{A}^\alpha_{L,K}$ exactly when $T(f) = 0$. We will use the following lemma:

\begin{lemma}
    \label{lem:2.1.3}
    If $\pi_K \mid f(x)$ in $\mathcal{O}_K[[x]]$, then $\pi_K^2 \mid T(f)$.
\end{lemma}

\noindent
Proof: we are working under the assumption that $\pi_L \mid \alpha$. It follows that $\pi_K^2 \mid [\alpha]_L(f([\pi_K](x)))$ if $\pi_K \mid f(x)$. Furthermore, assuming $\pi_K \mid f(x)$ and noting \hyperref[lem:2.1.2]{Lemma 2.6} implies we have
$$\sum^{\text{LT},L}_{z \in \textfrak{F}_0(F_K)} f(x \oplus_K z) \equiv \sum_{z\in\textfrak{F}_0(F_K)} f(x\oplus_K z) \equiv \mathscr{L}_{F_K}(f)([\pi_K](x)) \mod \pi_K^{q_L}$$

\noindent
Then because $q_L \ge 2$ we can conclude $T(f) \equiv 0 \mod \pi_K^2$ whenever $\pi_K \mid f(x)$ by Lemma 6 of \cite{Coleman1}.
\\
\\
We have $T(r) = 0$ because $r(x) \in \mathscr{A}^\alpha_{L,K}$. We also have
$$T(s) = T(s) \ominus_L T(r) = T(s \ominus_L r)$$

\noindent
We know $s(x) \ominus_L r(x) \equiv 0 \mod \pi_K$. We conclude that $T(s) \equiv 0 \mod \pi_K^2$ from \hyperref[lem:2.1.3]{Lemma 2.7}. This implies one direction of \hyperref[thm:2.1.1]{Theorem 2.5}.
\\
\\
To prove the converse of the above we take some series $s(x) \in \mathcal{O}_K[[x]]$ whose coefficients are all units. We suppose
$$\sum^{\text{LT},L}_{z \in \textfrak{F}_0(F_K)} s(x \oplus_K z) \equiv [\alpha]_L(s([\pi_K]_K(x)) \mod \pi_K^2$$

\noindent
and we use this to show there exists some $r(x) \in \mathscr{A}^\alpha_{L,K}$ such that $r(x) \equiv s(x) \mod \pi_K$. To prove this it suffices to find some $t(x) \in \pi_K\mathcal{O}_K[[x]]$ such that $T(s\ominus_L t) = 0$. Consider $T(t(x))$ for $t(x) \in \pi_K\mathcal{O}_K[[x]]$. As long as $q_L \ge 3$ and $\pi_K^2 \mid \pi_L$ we have
$$[\alpha]_L(t([\pi_K]_K(x))) \equiv 0 \mod \pi_K^3$$

\noindent
It follows that
$$T(t) \equiv \sum^{\text{LT},L}_{z \in \textfrak{F}_0(F_K)} t(x \oplus_K z) \mod \pi_K^3$$

\noindent
Then by \hyperref[lem:2.1.2]{Lemma 2.6} we have
$$T(t) \equiv \sum^{\text{LT},L}_{z \in \textfrak{F}_0(F_K)} t(x \oplus_K z) \equiv \sum_{z \in \textfrak{F}_0(F_K)} t(x \oplus_K z) \mod \pi_K^3$$

\noindent
assuming $q_L \ge 3$. Our goal is to pick a series $t_1(x) \in \pi_K\mathcal{O}_K[[x]]$ such that
$$T(t_1) \equiv \sum_{z \in \textfrak{F}_0(F_K)} t_1(x \oplus_K z) \equiv T(s) \mod \pi_K^3$$

\noindent
To show that it is possible to pick such a $t_1$ note that $T(s) = s_1([\pi_K]_K(x))$ for some series $s_1(x) \in \pi_K^2\mathcal{O}_K[[x]]$ by Lemma 3 of \cite{Coleman1}. The map $\mathscr{L}_{F_K} : \mathcal{O}_K[[x]] \rightarrow \pi_K\mathcal{O}_K[[x]]$ is surjective by the proof of Lemma 16 in \cite{Coleman2}. Therefore we can choose $t_1(x)$ to be any series in $\pi_K\mathcal{O}_K[[x]]$ satisfying $\mathscr{L}_{F_K}(t_1) = s_1$. For such a choice of $t_1$ we have
$$T(t_1) \equiv \sum_{z \in \textfrak{F}_0(F_K)} t_1(x \oplus_K z) \equiv \mathscr{L}_{F_K}(t_1)([\pi_K]_K(x)) \equiv s_1([\pi_K]_K(x)) \mod \pi_K^3$$

\noindent
and we conclude that for this choice of $t_1$ we have
$$T(t_1) \equiv T(s) \mod \pi_K^3$$

\noindent
Now fix some $N \ge 1$. Suppose for each $n$ with $1 \le n \le N$ we have $t_n \in \pi_K^n\mathcal{O}_K[[x]]$ and suppose that
$$T(s \ominus_L (t_1 \oplus_L \ldots \oplus_L t_n)) \equiv 0 \mod \pi_K^{n+2}$$

\noindent
We show we can pick a series $t_{N+1}(x) \in \pi_K^{N+1}\mathcal{O}_K[[x]]$ such that
$$T(s \ominus_L (t_1 \oplus_L \ldots \oplus_L t_{N+1})) \equiv 0 \mod \pi_K^{N+3}$$

\noindent
First note we have
$$T(s \ominus_L (t_1 \oplus_L \ldots \oplus_L t_N)) \equiv 0 \mod \pi_K^{N+2}$$

\noindent
If $t_{N+1}(x)$ is any series in $\pi_K^{N+1}\mathcal{O}_K[[x]]$ we have
$$T(t_{N+1}) \equiv \sum^{\text{LT},L}_{z \in \textfrak{F}_0(F_K)} t_{N+1}(x \oplus_K z) \equiv \mathscr{L}_{F_K}(t_{N+1})([\pi_K]_K(x)) \mod \pi_K^{N+3} $$

\noindent
Again by Lemma 3 of \cite{Coleman1} we have there exists some series $s_N(x) \in \pi_K^{N+2}\mathcal{O}_K[[x]]$ such that
$$T(s \ominus_L (t_1 \oplus_L \ldots \oplus_L t_N)) = s_N([\pi_K]_K(x))$$

\noindent
It follows that if we pick $t_{N+1}(x)$ to be any series in $\pi_K^{N+1}\mathcal{O}_K[[x]]$ such that $\mathscr{L}_{F_K}(t_{N+1}) = s_N$ we get that
$$T(t_{N+1}) \equiv \mathscr{L}_{F_K}(t_{N+1})([\pi_K]_K(x)) \equiv s_N([\pi_K]_K(x)) \equiv T(s \ominus_L (t_1 \oplus_L \ldots \oplus_L t_N)) \mod \pi_K^{N+3}$$

\noindent
For this choice of $t_{N+1}(x)$ we conclude that
$$T(s \ominus_L (t_1 \oplus_L \ldots \oplus_L t_{N+1})) \equiv 0 \mod \pi_K^{N+3}$$

\noindent
as desired. The sequence of partial sums $s \ominus_L (t_1 \oplus_L \ldots \oplus_L t_N)$ converges coefficientwise to some series $r(x)$. $T(r)$ must be zero because the expressions $T(s \ominus_L (t_1 \oplus_L \ldots \oplus_L t_N))$ converge coefficientwise to zero. We conclude for this series $r(x)$ we have $r(x) \in \mathscr{A}^\alpha_{L,K}$ and $r(x) \equiv s(x) \mod \pi_K$. This completes the proof of \hyperref[thm:2.1.1]{Theorem 2.5}.
\\
\\
We now prove the following lemma which allows us to construct some series $r(x) \in \mathscr{A}^\alpha_{L,K}$ where $\alpha = q_K/\pi_L^{f_{K/L}}$ such that $r(x)$ does not vanish mod $\pi_K$.

\begin{lemma}
    \label{lem:2.1.4}
    Suppose $s(x) \in \mathcal{O}_K[[x]]$ is such that all coefficients of $s$ are units. We also assume the extension $K$ over $L$ is ramified, so $|\pi_L| < |\pi_K|$. Then if there exists some $s_0(x) \in \mathcal{O}_K[[x]]$ such that $s(x) = s_0(x^{q_K})$ we have there exists some $r(x) \in \mathscr{A}^\alpha_{L,K}$ such that $r(x) \equiv s(x) \mod \pi_K$.
\end{lemma}

\noindent
Proof: for $s(x) = s_0(x^{q_K})$ one can check that 
$$\sum^{\text{LT},L}_{z \in \textfrak{F}_0(F_K)} s(x \oplus_K z) \equiv [\alpha]_L(s([\pi_K]_K(x)) \mod \pi_K^2$$

\noindent
We then conclude the lemma by applying \hyperref[thm:2.1.1]{Theorem 2.5}. We check for any nonzero $z \in \textfrak{F}_0(F_K)$ we have the following congruence whenever $\pi_K^2 \mid p$:
$$(x \oplus_K z)^{q_K} \equiv x^{q_K} \mod z^{q_K}$$

\noindent
The above congruence holds because $y$ divides every nonlinear term of $F(x,y)$ and because $\pi_K^2 \mid {q_K \choose k}$ for every index $0 < k < q_K$ under the assumption that $\pi_K^2 \mid p$. We conclude from the above that
$$\sum^{\text{LT},L}_{z \in \textfrak{F}_0(F_K)} s(x \oplus_K z) \equiv \sum^{\text{LT},L}_{z \in \textfrak{F}_0(F_K)} s_0((x \oplus_K z)^{q_K}) \equiv \sum^{\text{LT},L}_{z \in \textfrak{F}_0(F_K)} s_0(x^{q_K}) \mod z^{q_K}$$

\noindent
Since the series $\sum^{\text{LT},L}_{z \in \textfrak{F}_0(F_K)} s(x \oplus_K z)$ lives in $\mathcal{O}_K[[x]]$ we conclude that the previous congruence is actually a congruence mod $\pi_K^2$. We conclude that
$$\sum^{\text{LT},L}_{z \in \textfrak{F}_0(F_K)} s(x \oplus_K z) \equiv \sum^{\text{LT},L}_{z \in \textfrak{F}_0(F_K)} s_0(x^{q_K}) \equiv [q_K]_L(s(x)) \mod \pi_K^2$$

\noindent
We also have
$$[\alpha]_L(s([\pi_K]_K(x))) = [\alpha]_L(s_0([\pi_K]_K(x)^{q_K}))$$

\noindent
and $[\pi_K]_K(x)^{q_K} \equiv x^{q_K^2} \mod \pi_K^2$. We conclude
$$[\alpha]_L(s([\pi_K]_K(x))) \equiv [\alpha]_L(s_0(x^{q_K^2})) \equiv [\alpha]_L(s(x^{q_K})) \mod \pi_K^2$$

\noindent
Then note that $[\pi_L^{f_{K/L}}]_L(x) \equiv x^{q_K} \mod \pi_K^2$ assuming $|\pi_L| < |\pi_K|$. It follows that
$$[\alpha]_L(s([\pi_K]_K(x))) \equiv [\alpha]_L([\pi_L^{f_{K/L}}](s(x)) \equiv [q_K]_L(s(x)) \mod \pi_K^2$$

\noindent
We conclude that if $s(x) = s_0(x^{q_K})$ for some $s_0(x) \in \mathcal{O}_K[[x]]$ then $s(x)$ satisfies the conditions of \hyperref[thm:2.1.1]{Theorem 2.5}. It follows that for such $s(x)$ we have there exists $r(x) \in \mathscr{A}^\alpha_{L,K}$ such that $r(x) \equiv s(x) \mod \pi_K$, and this concludes the proof of \hyperref[lem:2.1.4]{Lemma 2.8}.
\\
\\
Together \hyperref[lem:2.1.4]{Lemma 2.8} and that $\mathscr{A}^\alpha_{L,K,\pi_L}$ is isomorphic to the $\alpha$ eigenspace of Coleman's trace operator in $\pi_L\mathcal{O}_K[[x]]$ provide a partial answer to the problem of classifying all series in $\mathscr{A}^\alpha_{L,K}$ when $\alpha = q_K/\pi_L^{f_{K/L}}$. It remains to classify all series in the module $\mathscr{A}^\alpha_{L,K}$ when $\alpha = q_K/\pi_L^{f_{K/L}}$.

\noindent

\section{Some additional eigenspaces of Coleman's trace operator}
\label{sec:section3}

\noindent
Let $\lambda \in \mathcal{O}_K$ be such that $\pi_K \mid \lambda$. In \cite{DiCapuaKolyvagin} the $\lambda$-eigenspace of $\mathscr{L}_{F_K}$ was shown to be isomorphic to the kernel of $\mathscr{L}_{F_K}$ in $\mathcal{O}_K[[x]]$. The goal of this section is to show the the same result holds when we replace the constant $\lambda$ with a power series $\alpha(x) \in \pi^2_K\mathcal{O}_K[[x]]$.
\\
\\
For a given $\alpha(x) \in \pi_K\mathcal{O}_K[[x]]$ we let $\mathscr{E}^\alpha_K$ denote the $\mathcal{O}_K$-module of all series $f(x) \in \mathcal{O}_K[[x]]$ such that
$$\mathscr{L}_{F_K}(f) = \alpha(x)f(x)$$

\noindent
Let $\mathscr{C}_K$ denote the kernel of Coleman's trace operator in $\mathcal{O}_K[[x]]$ as an $\mathcal{O}_K$-module. We define the map $\rho_\alpha : \mathscr{E}^\alpha_K \rightarrow \mathscr{C}_K$ by
$$\rho_\alpha(f) = f(x) - \frac{k(x)\alpha([\pi_K]_K(x))f([\pi_K]_K(x))}{\pi_Kw([\pi_K]_K(x))}$$

\noindent
We are able to prove the following:
\begin{theorem}
    \label{thm:3.0.1}
    The map $\rho_\alpha : \mathscr{E}^\alpha_K \rightarrow \mathscr{C}_K$ is an isomorphism of $\mathcal{O}_K$-modules.
\end{theorem}

\noindent
In the above definition the series $k(x)$ and $w(x)$ are the same $k$ and $w$ from the previous section. The proof is mostly the same as when $\alpha(x)$ is some constant divisible by $\pi_K$.
\\
\\
It is clear that for $f(x) \in \mathcal{O}_K[[x]]$ the expression $\rho_\alpha(f) \in \mathcal{O}_K[[x]]$ because $\pi_K \mid \alpha(x)$ and because $w([\pi_K]_K(x))$ is a unit in the ring $\mathcal{O}_K[[x]]$. We must check $\mathscr{L}_{F_K}(\rho_\alpha(f)) = 0$ in order to show $\rho_\alpha(f) \in \mathscr{C}_K$ when $f \in \mathscr{E}^\alpha_K$.
\\
\\
Consider the expression 
$$\mathscr{L}_{F_K}(\rho_\alpha(f)) = \mathscr{L}_{F_K}(f) - \mathscr{L}_{F_K}(\frac{k(x) \alpha([\pi_K]_K(x)f([\pi_K]_K(x))}{\pi_Kw([\pi_K]_K(x))})$$

\noindent
We have $\mathscr{L}_{F_K}(f) = \alpha(x)f(x)$ because $f(x) \in \mathscr{E}^\alpha_K$. By the construction of the series $k$ and $w$ and the arguments in \cite{Coleman2} we know that
$$\mathscr{L}_{F_K}(\frac{k(x)\alpha([\pi_K]_K(x))f([\pi_K]_K(x))}{\pi_Kw([\pi_K]_K(x))}) = \alpha(x)f(x)$$

\noindent
Combining this with the previous expression for $\mathscr{L}_{F_K}(\rho_\alpha(f))$ gives us
$$\mathscr{L}_{F_K}(\rho_\alpha(f)) = \alpha(x)f(x) - \alpha(x)f(x) = 0$$

\noindent
so we conclude $\rho_\alpha(f) \in \mathscr{C}_K$. Next we would like to show the map $\rho_\alpha$ is injective.
\\
\\
In order to show $\rho_\alpha$ is injective we need to further stipulate that $\pi_K^2 \mid \alpha(x)$. Under these conditions, suppose we have some $f(x) \in \mathcal{O}_K[[x]]$ such that $\rho_\alpha(f) = 0$. If $f$ is nonzero we have there exists some nonnegative integer $N$ such that $\pi_K^N \mid f(x)$ and $\pi_K^{N+1}$ does not divide $f(x)$. Expanding out the definitin of $\rho_\alpha(f)$ in the equation $\rho_\alpha(f) = 0$ gives us
$$f(x) = \frac{k(x)\alpha([\pi_K]_K(x)) f([\pi_K]_K(x))}{\pi_Kw([\pi_K]_K(x))}$$

\noindent
Then because $\pi_K^2 \mid \alpha(x)$ one can check that $\pi_K^{N+1}$ divides the right side of the above equation. This is a contradiction because $\pi_K^{N+1}$ does not divide $f(x)$. We conclude that $\rho_\alpha$ is injective under the condition $\pi_K^2 \mid \alpha(x)$.
\\
\\
We now prove $\rho_\alpha : \mathscr{E}^\alpha_K \rightarrow \mathscr{C}_K$ is surjective. Given $h(x) \in \mathscr{C}_K$ we construct a series $f(x)$ such that $\rho_\alpha(f) = h$. We check that $f(x) \in \mathscr{E}^\alpha_K$.
\\
\\
We define a sequence of series $f_i(x) \in \mathcal{O}_K[[x]]$ such that the infinite sum $\sum_{i=1}^\infty f_i(x)$ converges coefficientwise to a series $f(x) \in \mathcal{O}_K[[x]]$. We construct the sequence $(f_i)$ in such a way that the partial sums $\rho_\alpha(\sum_{i=1}^Nf_i(x))$ converge coefficientwise to $h(x)$. Take $f_1(x) = h(x)$, so we have
$$\rho_\alpha(f_1) = h(x) - \frac{k(x) \alpha([\pi_K]_K(x))h([\pi_K]_K(x))}{\pi_Kw([\pi_K](x))}$$

\noindent
We take $f_2(x)$ to be $h(x) - \rho_\alpha(f_1)$, so that
$$f_2(x) = \frac{k(x) \alpha([\pi_K]_K(x)) h([\pi_K]_K(x))}{\pi_K w([\pi_K]_K(x))}$$

\noindent
If $M_h$ is the largest integer such that $\pi_K^{M_h} \mid h(x)$ we note that $\pi_K^{M_h + 1} \mid f_2(x)$. We also note that
$$\rho_\alpha(f_1 + f_2) = \rho_\alpha(f_1) + \rho_\alpha(f_2) = h(x) - \frac{k(x)\alpha([\pi_K]_K(x))f_2([\pi_k]_K(x))}{\pi_Kw([\pi_K]_K(x))}$$

\noindent
implying that $\pi_K^{M_h + 2} \mid h(x) - \rho_\alpha(f_1 + f_2)$.
\\
\\
Take some integer $N \ge 2$ and suppose for all $1 \le n \le N$ we have found series $f_n$ such that $\pi_K^{M_h + n - 1} \mid f_n(x)$ and such that
$$\pi_K^{M_h + n} \mid h(x) - \rho_\alpha(\sum_{i=1}^nf_i(x))$$

\noindent
We show it is possible to pick $f_{N+1}(x)$ such that $\pi_K^{M_h+N} \mid f_{N+1}(x)$ and
$$\pi_K^{M_h + N + 1} \mid h(x) - \rho_\alpha(\sum_{i=1}^{N+1}f_i(x))$$

\noindent
Take $f_{N+1}(x) = h(x) - \rho_\alpha(\sum_{i=1}^Nf_i(x))$. For this choice of $f_{N+1}(x)$ we get that $\pi_K^{M_h + N} \mid f_{N+1}(x)$. We also have
$$\rho_\alpha(f_{N+1}) = h(x) - \rho_\alpha(\sum_{i=1}^Nf_i(x)) - \frac{k(x)\alpha([\pi_K]_K(x)) f_{N+1}([\pi_K]_K(x))}{\pi_Kw([\pi_K]_K(x))}$$

\noindent
implying $\pi_K^{M_h + N +1} \mid h(x) - \rho_\alpha(\sum_{i=1}^{N+1}f_i(x))$.
\\
\\
For the sequence $(f_i(x))$ constructed above we have $\pi_K^{M_h + i -1} \mid f_i(x)$ for each $i$ implying the partial sums $\sum_{i=1}^Nf_i(x)$ converge coefficientwise to some $f(x) \in \mathcal{O}_K[[x]]$. It follows that the expressions $\rho_\alpha(\sum_{i=1}^Nf_i(x))$ also converge coefficientwise to $\rho_\alpha(f)$. We have already seen these expressions converge coefficientwise to $h(x)$, so we conclude that $\rho_\alpha(f) = h$.
\\
\\
We need to check that if $h(x) \in \mathscr{C}_K$ and $f(x) \in \mathcal{O}_K[[x]]$ satisfies $\rho_\alpha(f) = h(x)$ then $f(x) \in \mathscr{E}^\alpha_K$. We have
$$f(x) - \frac{k(x)\alpha([\pi_K]_K(x)) f([\pi_K]_K(x))}{\pi_Kw([\pi_K]_K(x))} = h(x)$$

\noindent
Applying $\mathscr{L}_{F_K}$ to both sides of the above equation gives 
$$\mathscr{L}_{F_K}(f) - \mathscr{L}_{F_K}( \frac{k(x)\alpha([\pi_K]_K(x)) f([\pi_K]_K(x))}{\pi_Kw([\pi_K]_K(x))} ) = 0$$

\noindent
From the arguments in the proof of Lemma 16 in \cite{Coleman2} we have the above is equivalent to
$$\mathscr{L}_{F_K}(f) - \alpha(x)f(x) = 0$$

\noindent
implying $f(x) \in \mathscr{E}^\alpha_K$. This concludes the proof that the map $\rho_\alpha : \mathscr{E}^\alpha_K \rightarrow \mathscr{C}_K$ is an isomorphism of $\mathcal{O}_K$-modules.

\noindent

\section{An application of the eigenspaces to some dihedral $(\varphi,\Gamma)$-modules}
\label{sec:section4}
\noindent
In this section we consider the $(\varphi,\Gamma)$-modules constructed by Berger in \say{A p-adic family of dihedral $(\varphi,\Gamma)$-modules} \cite{Berger1}. We refer to this paper for some notation and the description of the representations. Let $\mathbb{Q}_{p^2}$ denote the unique unramified quadratic extension of $\mathbb{Q}_p$. Let $\chi_2 : \text{Gal}(\overline{\mathbb{Q}_p}/\mathbb{Q}_{p^2}) \rightarrow \mathbb{Z}_{p^2}^\times$ be the character attached to the Lubin--Tate module over $\mathbb{Z}_{p^2}$ for the uniformizer $p$ of $\mathbb{Q}_{p^2}$. Given $x \in \mathbb{Z}_{p^2}^\times$, there is a unique way to write $x$ as $x = w(x) \langle x \rangle$ where $w(x)^{p^2 - 1} = 1$ and $\langle x \rangle \in 1 + p\mathbb{Z}_{p^2}$. For each $g \in \text{Gal}(\overline{\mathbb{Q}_p}/\mathbb{Q}_{p^2})$ we are then able to define $w_2(g) = w(\chi_2(g))$ and $\langle g \rangle_2 = \langle \chi_2(g) \rangle$. One of the main results of \cite{Berger1} is to explicitly construct the $(\varphi,\Gamma)$-modules associated to the representations $\text{Ind}^{\mathbb{Q}_p}_{\mathbb{Q}_{p^2}}(w_2(\cdot)^r\langle \cdot \rangle_2^s)$ where $r \in \mathbb{Z}$ and $s \in \mathbb{Z}_p$.
\\
\\
Still closely following \cite{Berger1} we give a definition of the Robba ring. Let $\mathcal{R}$ be the ring of all power series $f(x) = \sum_{n\in\mathbb{Z}} a_nx^n$ such that each $a_n \in \mathbb{Q}_p$ and $f(x)$ converges on some annulus $r_f \le |x| < 1$ where $r_f$ depends on $f$. We fix the following Frobenius action $\varphi$ on $\mathcal{R}$: $\varphi(f)(x) = f((1+x)^p-1)$.  
\\
\\
A $(\varphi,\Gamma)$-module over $\mathcal{R}$ is a free module of finite rank $d$ over $\mathcal{R}$ endowed with a semilinear Frobenius $\varphi$ such that the matrix of $\varphi$ lives in $\text{GL}_d(\mathcal{R})$ and a semilinear action of $\Gamma$ commuting with $\varphi$. By work of Fontaine, Cherbonnier-Colmez and Kedlaya there is an equivalence of categories between $p$-adic representations of $\text{Gal}(\overline{\mathbb{Q}_p}/\mathbb{Q}_p)$ and \'etale $(\varphi,\Gamma)$-modules over $\mathcal{R}$ \cite{Berger1}.
\\
\\
Let $Q(x) = \varphi(x)/x = ((1+x)^p-1)/x$. We also let $Q_2(x) = \varphi(Q(x))$. We will need Theorem A from \cite{Berger1}:

\begin{theorem*}
    (Berger 2011) If $r \in \mathbb{Z}$, $u \in \mathbb{Z}_p$, and $s = r + (p+1)u$, then the $(\varphi,\Gamma)$-module attached to $\text{Ind}^{\mathbb{Q}_p}_{\mathbb{Q}_{p^2}}(w_2(\cdot)^r \langle \cdot \rangle^s)$ has a basis in which the matrix for $\varphi$ is given by 

    $$\text{\emph{Mat}}(\varphi) = \begin{pmatrix}
        0 & 1\\
        Q_2^r(Q_2/Q^p)^u & 0
    \end{pmatrix}
    $$

    \noindent
    and the matrices for $\gamma \in \Gamma$ are given by

    $$\text{\emph{Mat}}(\gamma) = \begin{pmatrix}
        \chi(\gamma)^r\langle \chi(\gamma) \rangle^u(1 + O(x)) & 0\\
        0 & \chi(\gamma)^r\langle \chi(\gamma) \rangle^u(1 + O(x))
    \end{pmatrix} $$

    \noindent
    where the $1 + O(x)$ are two power series belonging to $1 + x\mathbb{Z}_p[[x]]$.
\end{theorem*}

\noindent
At this point we refer to \cite{FourquauxXie, BergerFourquaux} for how the $\psi$ operator of $p$-adic Hodge theory acts on a $(\varphi,\Gamma)$-module. First we define how $\psi$ acts on $\mathcal{R}$. Let $\textfrak{F}_0$ denote the set of roots of $\varphi(x) = (1+x)^p-1$. Then we define $\psi$ to be the unique $\mathbb{Q}_p$-linear operator satisfying

$$\psi(f)((1+x)^p-1) = \frac{1}{p}\sum_{z\in\textfrak{F}_0}f(x+z+xz)$$

\noindent
for all $f \in \mathcal{R}$. This $\psi$ is exactly \say{Coleman's trace operator} from \cite{Coleman1} scaled by a factor of $\frac{1}{p}$.
Now we can extend the $\psi$ operator on $\mathcal{R}$ to the $\psi$ operator on an arbitrary \'etale $(\varphi,\Gamma)$-module $M$ over $\mathcal{R}$. $\text{Mat}(\varphi)$ is invertible assuming $M$ is \'etale. Let $\{e_1, e_2, \ldots e_d\}$ be a basis for $M$. Then $\{\varphi(e_1),\varphi(e_2), \ldots, \varphi(e_d)\}$ is also a basis for $M$. For arbitrary $m \in M$ one can write $m = \sum_ir_i\varphi(e_i)$ where each $r_i \in \mathcal{R}$. We then define $\psi(m)$ to be $\psi(m) = \sum_i\psi(r_i)e_i$. Then this $\psi$ is the unique operator on $M$ such that $\psi(a\varphi(m)) = \psi(a)m$ and $\psi(\varphi(a)m) = a\psi(m)$ for all $a \in R$ and for all $m \in M$. $\psi$ commutes with $\Gamma$.
\\
\\
It is of particular interest to study the fixed points of $\psi$ on a given $(\varphi,\Gamma)$-module $M$. Let $M$ be the $(\varphi,\Gamma)$-module attached to $\text{Ind}^{\mathbb{Q}_p}_{\mathbb{Q}_{p^2}}(w_2(\cdot)^r\langle \cdot \rangle^s)$ from the previous theorem by Berger. Fix a basis $\{e_1,e_2\}$ of $M$ such that the matrix for $\varphi$ with respect to this basis is $\text{Mat}(\varphi) = \begin{pmatrix}
    0 & 1\\
    Q_2^r(Q_2/Q^p)^u & 0
\end{pmatrix}$. Then one can write $m = a\varphi(e_1) + b\varphi(e_2)$ where $a,b \in \mathcal{R}$. If we assume $\psi(m) = m$ then $\psi(m) = \psi(a)e_1 + \psi(b)e_2 = m = a\varphi(e_1) + b\varphi(e_2)$. By the work of Berger we know that

$$\varphi(e_1) = Q_2^r(Q_2/Q^p)^u e_2$$

\noindent
and

$$\varphi(e_2) = e_1$$
\noindent
It follows that

$$\psi(a)e_1 + \psi(b)e_2 = aQ_2^r(Q_2/Q^p)^ue_2 + be_1$$
\noindent
which can only happen if $\psi(a) = b$ and $\psi(b) = aQ_2^r(Q_2/Q^p)^u$. This of course implies $\psi^2(a) = Q^r_2(Q_2/Q^p)^ua$, so $a$ lives in the $Q_2^r(Q_2/Q^p)^u$-eigenspace of the operator $\psi^2$, and we clearly get a fixed point of $\psi$ in $M$ for every such series $a$ in the eigenspace. At this point we restrict our attention to $a$ with integral coefficients. 
\\
\\
Our main goals for the remainder of this section are the following: 
\phantomsection
\label{listofgoals}
\begin{enumerate}
    \item whenever $u \in \mathbb{Z}$ and $r \le -u(p+1)$ we are able to construct all fixed points of $\psi$ on $M$ having coefficients in the integral part of $\mathcal{R}$ with respect to the basis already fixed.
    \item whenever $u \in \mathbb{Z}$ and the conditions of (\hyperref[listofgoals]{1}) do not hold we are able to construct a family of fixed points of $\psi$ on $M$ having coefficients in the integral part of $\mathcal{R}$ with respect to the basis already fixed.
\end{enumerate}

\noindent
Let $\mathcal{R}^\text{int}$ denote the set of all $f \in \mathcal{R}$ such that all coefficients of $f$ live in $\mathbb{Z}_p$. For a series $\alpha(x) \in \mathcal{R}^\text{int}$ we let $\mathscr{E}^{\alpha,2}$ denote the $\alpha$ eigenspace of the operator $\psi^2$ in $\mathcal{R}^\text{int}$. Hence $\mathscr{E}^{\alpha,2} = \{ f \in \mathcal{R}^\text{int} \mid \psi^2(f)(x) = \alpha(x)f(x)\}$. For each $\alpha$ we define maps $\rho^{\alpha,2} : \mathscr{E}^{\alpha,2} \rightarrow \mathscr{C}_{\psi^2}(\mathcal{R}^\text{int})$, where $\mathscr{C}_{\psi^2}(\mathcal{R}^\text{int})$ denotes the kernel of $\psi^2$ in $\mathcal{R}^\text{int}$. We will just denote $\mathscr{C}_{\psi^2}(\mathcal{R}^\text{int})$ by $\mathscr{C}_{\psi^2}$ when there is no risk of confusion.
\\
\\
We will need a lemma about the structure of $\mathscr{C}_{\psi^2}$ as a $\mathbb{Z}_p$-module.

\begin{lemma}
    Let $\mathscr{C}(\mathcal{R}^\text{int}) = \mathscr{C}$ denote the kernel of $\psi$ in $\mathcal{R}^\text{int}$. Then $\mathscr{C}_{\psi^2}$ and $\mathscr{C}\times\mathscr{C}$ are isomorphic as $\mathbb{Z}_p$-modules.
\end{lemma}

\noindent
Proof: let $k(x)$ and $w(x)$ be the series from \hyperref[sec:section2]{section 2} with the same names, specialized for the multiplicative formal group law. Define a map $t : \mathscr{C} \times \mathscr{C} \rightarrow \mathscr{C}_{\psi^2}$ by letting 

$$t(h,g) = h(x) + \frac{k(x)\varphi(g)(x)}{\varphi(w)(x)}$$

\noindent
Then 

$$\psi(t(h,g)) = \psi(h(x)) + \psi(\frac{k(x)\varphi(g)(x)}{\varphi(w)(x)}) = g(x)$$

\noindent
by Coleman's lemma in \cite{Coleman2}. It follows we have $\psi^2(t(h,g)) = 0$ for all choices of $(h,g)\in \mathscr{C}\times\mathscr{C}$. We need to show that $t$ is injective and surjective. To show $t$ is injective note that if $t(h,g) = 0$ then by the above computation $g(x) = 0$. It follows that we must have $h(x) = 0$ as well. This shows $t$ is injective. 
\\
\\
To prove $t$ is surjective, take an arbitrary $l(x) \in \mathscr{C}_{\psi^2}$. $\psi(l)(x) = l_2(x) \in \mathscr{C}$, so we can write 

$$l(x) = l_1(x) + \frac{k(x)\varphi(l_2)(x)}{\varphi(w)(x)} $$

\noindent
where $l_1(x)$ also lives in $\mathscr{C}$. This completes the proof that the map $t$ is surjective.
\\
\\
First we prove (\hyperref[listofgoals]{1}). Assume $u \in \mathbb{Z}$ and $r \le -u(p+1)$. Our goal is to find all series $a \in \mathcal{R}^\text{int}$ such that $\psi^2(a) = Q_2^r(Q_2/Q^p)^ua = Q_2^{r+u}Q^{-up}a$. We will need the following lemma:

\begin{lemma}
    \label{lem:4.0.2}
    Let $\alpha(x)$ be an arbitrary series in $\mathcal{R}^\text{int}$. Let $\beta(x) = Q_2(x)Q(x)\alpha(x)$. Then the explicit map $\rho_{\alpha,\beta} : \mathscr{E}^{\alpha,2} \rightarrow \mathscr{E}^{\beta,2}$ defined by $\rho_{\alpha,\beta}(f) = f(x)/\varphi^2(x)$ is an isomorphism of $\mathbb{Z}_p$-modules.
\end{lemma}

\noindent
Proof: let $f(x) \in \mathscr{E}^{\alpha,2}$. Then

$$\psi^2(f)(x) = \alpha(x)f(x)$$

\noindent
Now $\psi(f/\varphi^2(x)) = \frac{1}{\varphi(x)}\psi(f)$, so that

$$\psi^2(f/\varphi^2(x)) = \psi( \frac{1}{\varphi(x)} \psi(f)) = \frac{1}{x}\psi^2(f) = \frac{1}{x}\alpha(x)f(x) = \frac{\varphi^2(x)}{x}\alpha(x)(f(x)/\varphi^2(x))$$

\noindent
So we have $\rho_{\alpha,\beta}(f) = f(x)/\varphi^2(x)$ lives in $\mathscr{E}^{\beta,2}$. The map is injective, and if we define $\rho_{\beta,\alpha} : \mathscr{E}^{\beta,2} \rightarrow \mathscr{E}^{\alpha,2}$ by $\rho_{\beta,\alpha}(g) = \varphi^2(x)g(x)$, then it is the case that $\rho_{\alpha,\beta}$ and $\rho_{\beta,\alpha}$ are inverses. We get that both maps are surjections, so we can conclude that both maps are isomorphisms of $\mathbb{Z}_p$-modules.
\\
\\
We are proving (\hyperref[listofgoals]{1}), so we are in the case where $-r - u(p+1) \ge 0$. Let 

$$\gamma = \frac{Q_2^{r+u}Q^{-up}}{(Q_2Q)^{r+u}} = Q^{-u(p+1) -r} \in \mathbb{Z}_p[[x]]$$

\noindent
Define the map $\rho_{\gamma,2} : \mathscr{E}^{\gamma,2} \rightarrow \mathscr{C}_{\psi^2}$ by

$$\rho_{\gamma,2}(f) = f(x) - \frac{k(x) \varphi( \frac{k(x)\varphi(\gamma)(x)\varphi(f)(x)}{\varphi(w)(x)} )}{\varphi(w)(x)} = f(x) - \frac{k(x) \varphi(k)(x) \varphi^2(\gamma)(x) \varphi^2(f)(x)}{\varphi(w)(x)\varphi^2(w)(x)}$$

\noindent
We must check $\psi^2(\rho_{\gamma,2}(f)) = 0$. First

$$\psi(\rho_{\gamma,2}(f)) = \psi(f) - \psi\left(\frac{k(x) \varphi(k)(x) \varphi^2(\gamma)(x) \varphi^2(f)(x)}{\varphi(w)(x)\varphi^2(w)(x)}\right) = \psi(f) - \frac{k(x) \varphi(\gamma)(x)\varphi(f)(x)}{\varphi(w)(x)}$$

\noindent
by the proof of Lemma 16 from \cite{Coleman2}. Then

$$\psi^2(\rho_{\gamma,2}(f)) = \psi^2(f) - \psi \left( \frac{k(x) \varphi(\gamma)(x)\varphi(f)(x)}{\varphi(w)(x)} \right) = \psi^2(f) - \gamma(x)f(x) =0$$

\noindent
again by the proof of Lemma 16 from \cite{Coleman2}. Let $\mathscr{C}_{\psi^2}^+ = \mathscr{C}_{\psi^2} \bigcap \mathbb{Z}_p[[x]]$. We also prove that the image $\rho_{\gamma,2}(\mathscr{E}^{\gamma,2})$ contains every $h \in \mathscr{C}^+_{\psi^2}$.
\\
\\
We prove the following theorem:

\begin{theorem}
    \label{thm:4.0.3}
    Let $h(x)$ be an arbitrary power series in $\mathscr{C}^+_{\psi^2}$. Then there exists a unique $f_h \in \mathbb{Z}_p[[x]]$ such that $\rho_{\gamma,2}(f_h) = h$
\end{theorem}

\noindent
Proof: our goal is to construct a power series $f_h(x)$ such that

$$f_h(x) - \frac{k(x)\varphi(k)(x)\varphi^2(\gamma)(x)\varphi^2(f_h)(x)}{\varphi(w)(x) \varphi^2(w)(x)} = h(x)$$

\noindent
We accomplish this by constructing a sequence of series $(f_n)_{n\ge 0}$ with each $f_n(x) \in \mathbb{Z}_p[[x]]$ such that 

$$\rho_{\gamma,2}(\sum_{n=0}^Nf_n(x)) - h(x)$$

\noindent
is divisible by $x^{(N+1)2(p-1)}$ and such that each $f_n(x)$ is divisible by $x^{2n(p-1)}$. Once we have such a sequence 

$$\lim_{N\rightarrow\infty}\sum_{n=0}^Nf_n$$

\noindent
will converge coefficientwise to some series in $\mathbb{Z}_p[[x]]$, and we will be able to take $f_h$ to be this limit.
\\
\\
Take $f_0$ to be $f_0(x) = h(x)$. Then

\begin{align*}
    \rho_{\gamma,2}(f_0) - h(x) = f_0(x) - \frac{k(x)\varphi(k)(x)\varphi^2(\gamma)(x)\varphi^2(f_0)(x)}{\varphi(w)(x) \varphi^2(w)(x)} - h(x)\\
    =h(x) - \frac{k(x)\varphi(k)(x)\varphi^2(\gamma)(x)\varphi^2(h)(x)}{\varphi(w)(x) \varphi^2(w)(x)} - h(x)\\
    = -\frac{k(x)\varphi(k)(x)\varphi^2(\gamma)(x)\varphi^2(h)(x)}{\varphi(w)(x) \varphi^2(w)(x)}
\end{align*}

\noindent
Since $x^{p-1} \mid k(x)$, we get that $\rho_{\gamma,2}(f_0) - h(x)$ is divisible by $x^{2(p-1)}$. Now suppose we have found $f_n$ up to $n = N$ satisfying the conditions, so that we have

$$\rho_{\gamma,2}(\sum_{n=0}^Nf_n(x)) - h(x)$$

\noindent
is divisible by $x^{(N+1)2(p-1)}$ and for each $n$ we have $x^{2n(p-1)}\mid f_n(x)$. We can define $f_{N+1}$ to be

$$f_{N+1}(x) = h(x) - \rho_{\gamma,2}(\sum_{n=0}^Nf_n(x))$$

\noindent
For this choice of $f_{N+1}$ it is immediate that we have $x^{(N+1)2(p-1)}\mid f_{N+1}(x)$. Furthermore we have

\begin{align*}
\rho_{\gamma,2}(\sum_{n=0}^{N+1}f_{n}(x)) - h(x) = \rho_{\gamma,2}(f_{N+1}(x)) + \rho_{\gamma,2}(\sum_{n=0}^Nf_n(x)) - h(x)\\
= f_{N+1}(x) - \frac{k(x)\varphi(k)(x)\varphi^2(\gamma)(x)\varphi^2(f_{N+1})(x)}{\varphi(w)(x) \varphi^2(w)(x)} + \rho_{\gamma,2}(\sum_{n=0}^Nf_n(x)) - h(x)\\
= - \frac{k(x)\varphi(k)(x)\varphi^2(\gamma)(x)\varphi^2(f_{N+1})(x)}{\varphi(w)(x) \varphi^2(w)(x)}
\end{align*}

\noindent
which is divisible by $x^{(N+1)2(p-1) + 2(p-1)} = x^{(N+2)2(p-1)}$ because $f_{N+1}(x) $ is divisible by $x^{(N+1)2(p-1)}$ and hence so is $\varphi^2(f_{N+1}(x))$.
\\
\\
Now the sequence of partial sums

$$(\sum_{n=0}^Nf_n(x))$$

\noindent
must converge coefficientwise as $N \rightarrow \infty$ to some series which we denote by $f_h$. We have shown above that the sequence of series

$$\rho_{\gamma,2}(\sum^N_{n=0} f_n(x))$$

\noindent
converges coefficientwise to $h(x)$. It is clear that the same sequence must also converge coefficientwise to $\rho_{\gamma,2}(f_h)$ from the definition of $\rho_{\gamma,2}$, so we can conclude that $\rho_{\gamma,2}(f_h) = h(x)$. 
\\
\\
To show $f_h$ is unique we will show the map $\rho_{\gamma,2}$ is injective on $\mathbb{Z}_p[[x]]$. Suppose there is some nontrivial $f$ in the kernel so that

$$f(x) - \frac{k(x)\varphi(k)(x)\varphi^2(\gamma)(x)\varphi^2(f)(x)}{\varphi(w)(x)\varphi^2(w)(x)} = 0$$

\noindent
Note the above would imply that if $x^N$ is the largest power of $x$ dividing $f(x)$ then $x^{N+2(p-1)}$ also divides $f$, a contradiction. Therefore we conclude $\rho_{\gamma,2}$ is injective. The same map is clearly $\mathbb{Z}_p$-linear. \hyperref[thm:4.0.3]{Theorem 4.3} is proved.
\\
\\
Now if $\alpha = Q_2^{r+u}Q^{-up}$ we can iterate \hyperref[lem:4.0.2]{Lemma 4.2} on the series constructed in the $\gamma$-eigenspace of $\psi^2$ to get series in the $\alpha$-eigenspace of $\psi^2$. Furthermore if $r + u = -up$ so that $r = -u(p+1)$ then $\gamma = 1$, and there is already a classification of the $\gamma$-eigenspace of $\psi$ due to Colmez (the same argument can be slightly modified to find the $\gamma$-eigenspace of $\psi^2$ in this specific case).
\\
\\
We now move to classifying arbitrary $\mathscr{E}^{\gamma,2}$ given the conditions $u \in \mathbb{Z}$ and $r \le -u(p+1)$. We accomplish this by showing that if $f \in \mathscr{E}^{\gamma,2}$ then $f \in \mathbb{Z}_p[[x]]$. 
\\
\\
Let $f$ be in $\mathcal{R}^\text{int}$ such that

$$\mathscr{L}^2(f) = p^2\gamma(x)f(x)$$

\noindent
where $\gamma(x) = Q^{-u(p+1)-r}(x) = (((x+1)^p-1)/x)^{-u(p+1) - r}$. Fix a sequence $(u_n)$ where each $u_n$ is an $n$-th level torsion point of the multiplicative formal group law, $u_0 \neq 0$, and $(u_{n+1}+1)^p-1 = u_n$ for each $n$. Now $\gamma(x)$ has a multiplicative inverse which converges on the annulus $|u_0| < |x| < 1$. We have

$$\sum_{z\in\textfrak{F}_1}f(x\oplus z) = p^2\varphi^2(\gamma)(x)\varphi^2(f)(x)$$

\noindent
where $\textfrak{F}_1$ denotes the set $\{z \mid \varphi^2(x) \text{ at } z \text{ is } 0 \}$. $\varphi^2(\gamma)(x)^{-1}$ converges on the annulus $|u_2| < |x| < 1$, so if $f$ converges on $|u_n| < |x| < 1$ then $f$ converges on $|u_{n-2}| < |x| < 1$ for all $n \ge 2$. This is also because $\varphi(x)$ takes the annulus $|u_n| < |x| < 1$ to the annulus $|u_{n-1}| < |x| < 1$ for such indices $n$. This implies $f$ converges on $|u_0| < |x| < 1$. This also implies $\gamma(x)f(x)$ converges on $|p^2u_0|<|x|<1$.
\\
\\
At this point we know the series
{\abovedisplayskip=0pt\belowdisplayskip=0pt
\[
\sum_{z\in\textfrak{F}_1}f(x\oplus z) = p^2\varphi^2(\gamma)(x)\varphi^2(f)(x)
\]}
\noindent
converges on the annulus $|u_0| < |x| < 1$. In particular this annulus contains $u_1$, implying that $\gamma(x)f(x)$ is defined at zero since $\varphi^2(x)$ at $u_1$ is zero. Therefore we have $\gamma(x)f(x) \in \mathbb{Z}_p[[x]]$ and since $\gamma(x)$ has a formal power series inverse living in $K[[x]]$, we must have $f(x) \in \mathbb{Z}_p[[x]]$ as well. 
\\
\\
The above combined with \hyperref[thm:4.0.3]{Theorem 4.3} implies that the map $\rho_{\gamma,2}$ is an isomorphism of $\mathbb{Z}_p$-modules between $\mathscr{E}^{\gamma,2}$ and $\mathscr{C}^+_{\psi^2}$. This completes the classification of all series $a \in \mathcal{R}^\text{int}$ such that
{\abovedisplayskip=0pt\belowdisplayskip=0pt
\[
\psi^2(a) = Q_2^{r+u}Q^{-up}a
\]}
\noindent
under the conditions $u \in \mathbb{Z}$ and $r \le -u(p+1)$. We have proven the following theorem:

\begin{theorem}
    \label{thm:4.0.4}
    Let $r, u \in \mathbb{Z}$ and $s = r + (p+1)u$. Assume $s \le 0$. Then the $(\varphi,\Gamma)$-module attached to $\text{\emph{Ind}}^{\mathbb{Q}_p}_{\mathbb{Q}_{p^2}}(w_2(\cdot)^r \langle \cdot \rangle^s)$ has a basis for which the set of fixed points of the $\psi$-operator with integral coefficients with respect to this basis is isomorphic to $\mathscr{C}_{\psi^2}^+$ as a $\mathbb{Z}_p$-module. 
\end{theorem}

\noindent
Next we prove (\hyperref[listofgoals]{2}), so we are in the case where $u \in \mathbb{Z}$ and there are no longer any restrictions on $u(p+1)+r = s$. We define $\gamma$ exactly the same way as it is defined in the proof of (\hyperref[listofgoals]{1}), so that
\[
\gamma = \frac{Q_2^{r+u}Q^{-up}}{(Q_2Q)^{r+u}} = Q^{- u(p+1)-r} \in \mathcal{R}^{\text{int}}
\]

\noindent
Our first goal is to solve for a family of the series in $\mathscr{E}^{\gamma,2} = \{ f \in \mathcal{R}^{\text{int}} \mid \psi^2(f) = \gamma f \}$. We observe that if $f$ is such a series, then
\[
\sum_{z\in \textfrak{F}_1}f(x\oplus z) = p^2\varphi^2(\gamma)(x)\varphi^2(f)(x) = p^2\gamma(\varphi^2(x))f(\varphi^2(x))
\]

\noindent
Now $\gamma^{-1} \in \mathbb{Z}_p[[x]]$ (we can assume we are in the case $s > 0$ since we have already described the case $s\le 0$), so if $f$ converges on the annulus $|u_n| < |x| < 1$, then $f$ must also converge on the annulus $|u_{n-2}| < |x| < 1$ for all $n \ge 2$. This also follows because $\varphi$ takes the annulus $|u_n| < |x| < 1$ to the annulus $|u_{n-1}| < |x| < 1$. From this it follows that $f$ must converge on the annulus $|u_0| < |x| < 1$. Plugging in $u_1$ to the above equation implies $f$ must be defined at $x = 0$ (again because $\gamma^{-1} \in \mathbb{Z}_p[[x]]$). We get that all such $f \in \mathscr{E}^{\gamma,2}$ must live in $\mathbb{Z}_p[[x]]$.
\\
\\
We recall the map $\rho_{\gamma,2}: \mathscr{E}^{\gamma,2} \rightarrow \mathscr{C}_{\psi^2}$ defined earlier in this section:
\[
\rho_{\gamma,2}(f) = f(x) - \frac{k(x)\varphi(k)(x)\varphi^2(\gamma)(x)\varphi^2(f)(x)}{\varphi(w)(x)\varphi^2(w)(x)}
\]

\noindent
One initial observation is that, for such $f \in \mathscr{E}^{\gamma,2}$, one must have $\gamma^{-1}(x) \mid f(x)$ in $\mathbb{Z}_p[[x]]$. This is simply because $f\in \mathbb{Z}_p[[x]]$ and $\psi^2(f) = \gamma f$. We will be able to describe a submodule of the image $\rho_{\gamma,2}(\mathscr{E}^{\gamma,2}) \subseteq \mathscr{C}^+_{\psi^2}$. This is accomplished in the following:
\begin{theorem}
\label{theorem4.5}
Let $s = u(p+1)+r$. Then the image $\rho_{\gamma,2}(\mathscr{E}^{\gamma,2}) \subseteq \mathscr{C}^+_{\psi^2}$ contains each $h(x) \in \mathscr{C}^+_{\psi^2}$ such that $\varphi^2(x)^s \mid h(x)$ in $\mathbb{Z}_p[[x]]$.
\end{theorem}

\noindent
By \hyperref[lem:4.0.2]{lemma 4.2} we know that $\mathscr{E}^{\gamma,2}$ is isomorphic (as a $\mathbb{Z}_p$-module) to $\mathscr{E}^{\delta,2}$ where
\[
\delta(x) = Q_2^s(x) \in \mathbb{Z}_p[[x]]
\]

\noindent
We also know that this isomorphism is given by the map
\[
\rho_{\gamma,\delta} : \mathscr{E}^{\gamma,2} \rightarrow \mathscr{E}^{\delta,2} \text{ defined by } \rho_{\gamma,\delta}(f) \coloneq \frac{f(x)}{\varphi^2(x)^s}
\]

\noindent
It is possible to check that if $f\in\mathscr{E}^{\gamma,2}$ and $\rho_{\gamma,2}(f) = h_f$, then we must have 
\[
\rho_{\delta,2}(\rho_{\gamma,\delta}(f))(x) = \frac{h_f(x)}{\varphi^2(x)^s}
\]

\noindent
Similarly, if $g \in \mathscr{E}^{\delta,2}$ is arbitrary and $\rho_{\delta,2}(g) = h_g$, then
\[
\rho_{\gamma,2}(\rho_{\delta,\gamma}(g))(x) = \varphi^2(x)^s h_g(x)
\]

\noindent
where we use $\rho_{\delta,\gamma}$ to denote the inverse of $\rho_{\gamma,\delta}$. To prove \hyperref[theorem4.5]{Theorem 4.5} consider an arbitrary $h \in \mathscr{C}^+_{\psi^2}$. We check there is some $g_h \in \mathscr{E}^{\delta,2}$ such that $\rho_{\delta, 2}(g_h) = h$. Recall the definition of $\rho_{\delta,2}(g)$:
\[
\rho_{\delta,2}(g) = g(x) - \frac{k(x)\varphi(k)(x)\varphi^2(\delta)(x)\varphi^2(g)(x)}{\varphi(w)(x)\varphi^2(w)(x)}
\]

\noindent
Given arbitrary $h \in \mathscr{C}^+_{\psi^2}$, the proof that there exists some $g_h \in \mathscr{E}^{\delta,2}$ satisfying $\rho_{\delta,2}(g_h) = h$ works exactly the same as the proof of \hyperref[thm:4.0.3]{Theorem 4.3} since both $\delta \in \mathbb{Z}_p[[x]]$ and the choice of $\gamma$ used in \hyperref[thm:4.0.3]{Theorem 4.3} lives in $\mathbb{Z}_p[[x]]$. The same proof also shows that $g_h$ is unique.
\\
\\
Now once we have solved for $g_h \in \mathscr{E}^{\delta,2}$, we consider $\rho_{\delta, \gamma}(g_h) = \varphi^2(x)^s g_h \in \mathscr{E}^{\gamma,2}$. From our observations after the statement of \hyperref[theorem4.5]{Theorem 4.5}, we have that
\[
\rho_{\gamma,2}(\rho_{\delta,\gamma}(g_h))(x) = \varphi^2(x)^s h(x) \in \mathscr{C}^+_{\psi^2}
\]

\noindent
This completes the proof that the image $\rho_{\gamma,2}(\mathscr{E}^{\gamma,2})$ contains each $h(x) \in \mathscr{C}^+_{\psi^2}$ such that $\varphi^2(x)^s \mid h(x)$ in $\mathbb{Z}_p[[x]]$.

\subsection*{Acknowledgements}
I would like to thank Victor Kolyvagin for introducing me to this problem. I would like to thank Bryce Gollobit, Ryan Utke, Kioshi Morosin, and Gene Kopp for their detailed feedback on this paper. 

\subsection* {Funding}
I have no grant acknowledgements to declare.

\end{document}